\documentclass[titlepage]{article}[11pt]

\usepackage{latexsym, amssymb, amsmath,amsthm}
\usepackage[margin=1in]{geometry}
\usepackage{indentfirst}
\usepackage[font=footnotesize]{caption}
\usepackage{subcaption}
\usepackage[dvipsnames]{xcolor}
\usepackage{adjustbox}
\usepackage{lineno}
\usepackage[raggedrightboxes]{ragged2e}
\usepackage{comment}
\usepackage{soul}
\usepackage{longtable}
\usepackage{array}
\usepackage{url}
\usepackage{setspace}
\usepackage{apacite}
\usepackage{rotating}

\usepackage{multirow}

\title{Formulating human risk response in epidemic models: exogenous vs endogenous approaches}

\author{Leah LeJeune$^{1,2,a}$ \and Navid Ghaffarzadegan$^3$ \and Lauren M. Childs$^{1,2,b,*}$ \and Omar Saucedo$^{1,2,c,*,**}$}
\date{%
    $^1$Department of Mathematics, Virginia Tech\\
    225 Stanger Street Blacksburg, VA 24061-1026\\%
    $^a$leahlejeune@vt.edu\\
    $^b$lchilds@vt.edu\\
    $^c$osaucedo@vt.edu\\
    \vspace{1cm}
    $^2$Virginia Tech Center for the Mathematics of Biosystems, Virginia Tech\\
    765 West Campus Drive, Blacksburg, VA 24061\\
    \vspace{1cm}
    $^3$Department of Industrial and Systems Engineering, Virginia Tech\\
    7054 Haycock Rd, Falls Church, VA 22043\\
    navidg@vt.edu\\
    \vspace{1cm}
    *These authors contributed equally.\\
    **Corresponding author\\
[2ex]%
}

\begin{document}
\maketitle

\begin{abstract}
The recent pandemic emphasized the need to consider the role of human behavior in shaping epidemic dynamics. In particular, it is necessary to extend beyond the classical epidemiological structures to fully capture the interplay between the spread of disease and how people respond. Here, we focus on the challenge of incorporating change in human behavior in the form of ``risk response" into compartmental epidemiological models, where humans adapt their actions in response to their perceived risk of becoming infected. The review examines 37 papers containing over 40 compartmental models, categorizing them into two fundamentally distinct classes: exogenous and endogenous approaches to modeling risk response. While in exogenous approaches, human behavior is often included using different fixed parameter values for certain time periods, endogenous approaches seek for a mechanism internal to the model to explain changes in human behavior as a function of the state of disease. We further discuss two different formulations within endogenous models as implicit versus explicit representation of information diffusion. This analysis provides insights for modelers in selecting an appropriate framework for epidemic modeling.
\end{abstract}

\hspace{-2.35em}
\textbf{Keywords:} System dynamics, deterministic compartmental COVID-19 model, human behavior, risk response, epidemic modeling with endogenous formulations
\vspace{1em}
\hrule

\section{Introduction}

The role of human behavior in the spread of COVID-19 made abundantly clear the need to consider human response along with disease dynamics when simulating and predicting patterns of disease spread \cite{gentili2020challenges}. Without the inclusion of human response, model predictions for the extent and severity of COVID outbreaks were unreliable, specifically for long time horizons \shortcite{rahmandad2022enhancing,wirtz2021changing}. Human behavior impacted the pandemic in a multitude of ways, as many decisions were required on a daily basis. These include decisions such as whether to stay home within one's ``bubble" or go to work, to school, on errands, or to social events;  whether to mask when interacting with others; whether to test after a known exposure or after symptoms appear; whether to vaccinate; and so on. An additional set of decisions surrounding how to assess the state of the pandemic also existed: from which sources to obtain information, to whether to read or listen to the news, to whether to heed government recommendations, and more. Variation in human behavior changed the course of disease spread, which in turn caused changes in behavior. It is essential to capture these interconnected feedbacks in order to understand and predict pandemics.

As the pandemic unfolded and understanding grew about how disease spread was shaped by human behavior, new models emerged with a wide range of ways to incorporate human behavior \cite{sooknanan2023fomo,perra2021non}. These formulations were often simplistic, applying static choices to represent behavior for a period of time from a single region, such as the response of staying at home to a government imposed lockdown \cite{dick2021covid,childs2022modeling}. Models can be fit to and verified by past historical data, but this is specific to different scenarios and does not mean the model will give informative long-term predictions \cite{shankar2021systematic}. The lack of predictive ability of such models was noted as  motivation for the need to create  models that allow for responses that change in coordination with the state of the system \shortcite{rahmandad2021behavioral}, such as through a feedback loop incorporating human behavior with epidemiological modeling.

As an extension of a standard epidemiological model, \citeA{rahmandad2022enhancing} introduced the SEIR$b$ model, named for the standard susceptible (S), exposed (E), infectious (I), and removed (R) compartments, but also included a behavior ($b$) feedback. The focus of this model was understanding the influence of human behavior on disease spread via non-pharmaceutical interventions (NPIs). As opposed to the standard SEIR model without explicit behavioral links, it is able to capture several observed features of the pandemic such as multiple waves on an epidemic time scale and a shift towards an effective reproductive number of one. This incorporation of behavior also improved forecasting, with accuracy extending as far as 20 weeks ahead and, importantly, outperformed variations of the SEIR model, even when they included significantly more complexity \cite{rahmandad2022enhancing}.

One important application of forecasting is improved policy development to balance the cost of disease prevention with economic and social cost \shortcite{rahmandad2022missing}. For example, the SEIRb model challenged the hypothesis that focusing vaccination efforts on elderly populations is the most effective strategy for reducing deaths. Prioritizing those with low responsiveness rates (e.g., essential workers, high-contact individuals, communities with high death rates, or young communities with lower perceived risk) was suggested to be more effective than prioritizing high-risk individuals with high responsiveness rates such as elderly individuals \cite{rahmandad2022missing}, as was also observed from models that took careful consideration of contact structure \cite{childs2022modeling}. 

Thus, incorporating a human behavioral feedback loop into an existing epidemiological model structure is a promising way to advance disease forecasting and better inform policy creation to control disease spread while decreasing negative effects on society. Developing and exploring models of this type can help explain unusual disease patterns exhibited during the pandemic and is, thus, very relevant to prepare for future pandemics. 

When constructing a model incorporating human behavior, many decisions must be made to keep the model relevant but not overly complex. Human behavior has many facets that can range from spread of misinformation to vaccine hesitancy to adherence fatigue to name a few. 
These factors are important to consider when studying behavioral positive and negative feedback mechanisms as they play a significant role in the spread of infectious diseases. A positive feedback loop reinforces variable values, resulting in exponential growth. For example, word of mouth increases information spread in a social context. In contrast, a negative feedback loop regulates activities around a limit or equilibrium. For example, higher risks prompt more precaution, decreasing the risks \cite{sterman2000business}. The depth in which all of these behavioral feedback mechanisms can be studied is extensive. Thus, in this manuscript, we focus on ``risk response" and its incorporation into mathematical models. Risk response is driven by perceived risk, which is in turn driven by information related to the state of the disease. Perceived risk is the subjective assessment of the extent of risk associated with a specific event, which differs from actual risk in that humans feel the risk more than they understand the real risk, such as from assessing the current state of the disease analytically and objectively. This perception can cause people to adopt more NPI practices, which in turn decreases disease spread and thus forms a negative feedback loop. Conversely, such perception may make people less cautious, especially if cases are low. For example, before COVID-19 was widespread in the US, individuals were aware that future risk was approaching, but those in the proximity of the first infected individuals to appear in the US were not aware of their high risk until they were likely exposed. Thus, their perceived risk was lower than their actual risk, causing them to engage in more risky behaviors which contributed to significant disease spread. 

In this paper, we investigate how change in human behavior with respect to risk response has been integrated into mathematical models of COVID-19 spread and how this affects disease dynamics. Given the extensive literature on epidemic modeling post-pandemic, spanning hundreds of thousands of articles, we narrow our focus to a selective sample of 79 highly cited publications. This approach allows for a more in-depth analysis of the various modeling formulations. Our goal is to minimize our modeling efforts and use the existing literature as our data source, providing an article on how the scientific community models risk response. Our study is not intended to be a comprehensive report on incorporating human behavior. Instead, we specifically focus on compartmental models (formulated via ordinary differential equations) that consider change in human behavior throughout a pandemic, particularly concerning NPI adherence. By limiting our sample to more than 40 models compiled from 37 papers, we note and articulate two distinct classes of behavioral epidemic modeling in terms of risk response formulation as exogenous versus endogenous approaches. Furthermore, for the endogenous approach, we discuss varying representations of information diffusion. By constructing endogenous models using risk response mechanisms from the literature, we provide simulations which demonstrate the effects of risk response on disease dynamics, particularly the influence of the behavioral feedback loop in driving epidemic waves. Finally, we discuss our findings and contributions of our investigation, and consider future studies that go beyond the risk-response feedback loop.

\section{Categorizing risk response}

For clarity, Table \ref{table:definitions} lists the definitions central to this article. As previously stated, our focus is on formulating risk response as a behavioral construct that changes throughout a pandemic. In this paper, we particularly focus on mathematical representations of public response to rising or declining risks through change in compliance with NPIs. Such actions involve change in various protective practices such as masking, social distancing, quarantine, and isolation.

\begin{table}[ht!]
\small
\centering
\caption{\textbf{Glossary of Key Terms}} 
\begin{tabular}{|p{1.5in}|p{4.5in}|}
\hline
\textbf{Perceived risk} &  Individual and societal judgment about the likelihood of infection or harm from a disease, which may or may not be accurate.\\
\hline
\hline 
\textbf{Risk response} & Change in human behavior and compliance with NPIs as a result of change in perceived risk.\\
\hline
\hline
\textbf{Non-pharmaceutical interventions (NPIs)} & Policies implemented to slow the spread of disease which are non-medicinal, e.g. masking, social distancing, quarantining, isolating, but excluding vaccination and other drug-related treatments. \\
\hline
\hline
\textbf{Modeling risk response} &  Considering that higher (lower) perceived risk results in more (less) NPI compliance which in turn decreases (increases) infection rate.
\\
\hline
\hline
\textbf{Exogenous formulation of risk response} & Risk response is either assumed unchanging during the simulation, or changes based on predefined scenarios. It includes changing parameter values using external data sources, such as enacted policies or mobility data, or modelers' intuition with anecdotal evidence, such as the assumption of higher contact rate during holidays.  \\
\hline
\hline
\textbf{Endogenous formulation of risk response} &  Risk response is a function of the state of the disease, such as number of active cases or the subsequent deaths (as opposed to being a fixed parameter or a direct function of time). Consequently, not only does human behavior influence the spread of the disease, but the spread of the disease also influences human behavior, closing a feedback loop between disease and behavior.\\
\hline
\hline
\textbf{Implicit information diffusion} & Correlating NPI compliance with a disease state variable (or its lagged value). For example, in this approach we assume that average contact rate negatively correlates with lagged cases. The approach is agnostic to the mechanism by which individuals communicate such information among themselves.\\
\hline
\hline
\textbf{Explicit information diffusion} & Representing the spread of information about (or fear of) the epidemic through a population. In this approach, risk communication and the adoption of effective measures (such as masking) are modeled alongside the disease spread itself. 
\\
\hline

\end{tabular}
\label{table:definitions}
\end{table}

We differentiate formulations of human response to disease risk into two broad categories: exogenous and endogenous.  
Exogenous formulations incorporate human behavior by changes to fixed parameters, such that different scenarios are represented by different values. For the case of transmission rate (usually denoted as $\beta(\cdot)$, which may be a function of other variables), an exogenous formulation would be unaffected by the growth of the disease. In other words, change in the disease state does not drive individuals to alter their behavior. Thus, the transmission rate may change because the setting has changed, e.g. a government lockdown is imposed, but the state of the disease does not alter the value. As the choice of the parameter is dictated by external information, such as current government policy, it is not possible to forecast multiple waves without assuming changes to the parameter. Long-term forecasting is greatly hampered by the lack of feedback between disease level and human behavior. 

Endogenous formulations use a human behavior feedback loop to represent risk-driven response - that is, humans adjust their behavior after perceiving sufficient risk from information concerning the disease. For example, as the disease intensifies, people adhere more strictly to prevention policies, and as the disease abates, people relax their adherence, allowing for those quantities to increase again. Thus, the feedback loop between disease growth and human response drives disease waves in simulations, matching behavior observed during the pandemic.
These endogenous formulations can be further differentiated based on the type of information diffusion. Risk response is formulated using implicit information diffusion if a change in the state of the disease directly affects the disease transmission rate. This means that the transmission rate is a function of the disease state. 
In contrast, risk response is formulated using explicit information diffusion if there is a separate model or sub-model driving the behavior change. One formulation contains fully susceptible and less susceptible compartments (``split-susceptible"), with endogenous transfer between the compartments based on behavior changes resulting from awareness of the disease or in response to disease prevalence. Models with split-susceptibles often incorporate exogenous disease transmission. Another formulation is a coupled information/behavior-disease model where the information or behavior state variables are distinct from the disease model.

\section{Selection of sample of models}

Risk response is a complex human behavior and, thus, can be incorporated into epidemiological models in a myriad of ways. While our initial aim was a comprehensive description of how risk response has been included in disease models, given the vast amount of epidemic modeling papers, especially post COVID-19 pandemic (there are more than 5 million articles on just COVID-19 models according to Google Scholar), conducting a thorough systematic literature review was simply not feasible.
Instead, in order to gather different modeling formulations, we performed an extensive but focused review of literature to more efficiently assess different styles of incorporation of risk response in epidemiological models. 

\subsection{Focused literature search}

To that end, we started from one of the first literature reviews on epidemiological models that incorporate change in human behavior \cite{funk2010modelling}, which focused on conceptual differences in integrating human behavior into models. This source has been instrumental is shaping the direction of many relevant studies (e.g., \shortciteA{verelst2016behavioural}) and was very helpful in designing our study; however, our focus was on the mathematical mechanisms and formulations described there and in following works. As epidemic modeling has expanded enormously in the time since \citeA{funk2010modelling}, we also found a more recent review by \citeA{perra2021non}, which provided important insights on studies of incorporating change in human behavior during a disease outbreak. 

To capture a broad, albeit non-comprehensive and, thus manageable, swath of the literature, we used these two sources as an initial starting point. 
Specifically, we began with these two sources -- \cite{funk2010modelling} and \cite{perra2021non} -- and conducted a forward citation search in Google Scholar, restricting results from 2019-2023, which resulted in about 1,083 articles. However, despite citing the two source papers, focused on the need and methods for incorporation of human behavior in models, a majority of these articles are not primarily focused on modeling behavior change. Furthermore, many do not specifically include human behavior, particularly risk response, in their models. To narrow down search results to those which aligned with our focus (i.e., modeling risk response), we used the keywords of ``+math +model +risk +response +behavior +COVID". This narrowed the scope considerably and led to 168 articles, and we then excluded repeated search results. 
We further checked highly cited articles to assure a good coverage of our database, and as a part of this search, we found a few sources such as \shortciteA{verelst2016behavioural}, an important paper that builds on \citeA{funk2010modelling}, specifically mentioning exogenous approaches to behavior modeling.

Given the objective of this article - to analyze specific mathematical formulations - we then reviewed the articles to narrow down our sample size based on two major criteria. First, the models must  appear within peer-reviewed articles that were published since 2019, and, second, they must be  a mathematical model with an ``intentional” consideration of influence of human behavior, i.e. behaviors reflecting individual risk response. This is in contrast to an ``accidental” assumption of human behavior through mass action, incidence based exposures or well-mixed populations.  To assist with this evaluation, we conducted searches within articles using terms such as ``risk perception," ``risk response," ``change in behavior," and ``behavioral modeling." We read the papers carefully to make sure these criteria were satisfied, which brought the total to 79 papers. 

In addition to the forward citation search, we conducted a search on Google Scholar using several of those same terms as mentioned above such as ``risk perception," ``risk response," ``change in behavior," and ``behavioral modeling," and included several well-cited recent papers. Additional relevant papers  arising from various sources (e.g., recommendations from audiences at conferences or from academic peers) were included. This led to a total of 90 articles that included epidemic models. We then pared this collection down to 37 articles as described in the next section.

\subsection{Selection of models}

During the first stages of our analysis of the articles, we noted five main mathematical modeling approaches: compartmental models, agent-based/network models, statistical models, artificial intelligence models, and others including game theory. We noted that the first category (compartmental models) encompassed the largest category of the models - nearly half (41/90) - consistent with the arguments previously raised \cite{rahmandad2022enhancing}. In epidemiology, compartmental modeling is a convention that categorizes the population into distinct compartments based on disease status (e.g., in the classic SIR model into susceptible, infectious, and removed/recovered). While compartmental models frequently use differential equations to model the population flow between the compartments, these models can be mathematically formulated in a variety of ways.

We narrowed down our main focus to compartmental models (both deterministic and stochastic implementations), although it is important to note that these models were sometimes paired with other model types in an overarching model structure. We then carefully re-read the papers to assess precisely how response to risk is incorporated in the models. Since not all models use the terms ``risk response" in describing this mechanism, we specifically examined all features included in addition to the generic SIR/SEIR structure to determine how risk was mathematically formulated. In particular, models in which the rate of infection or contact rate was not fixed and changed for reasons other than herd immunity were considered as potential models for further analysis. This resulted in  a sample of over 40 mathematical models. This set of models, which by no means covers all potential modeling structures or approaches, highlights the major differences and individual contributions of various risk response formulations.   

In the remainder of this work, we first review the high-level differences between modeling approaches and provide a few emblematic examples. We then focus on the precise differences in mathematical formulations, including some illustrative simulation results. 

\section{Modeling approaches}\label{sec:formulationsriskresponse}

In this section, we report various approaches in the literature for formulating change in societal risk response, highlighting differences in model boundaries and information spread formulation.
 
\subsection{Model boundary and details}

One of the major differences across the models is related to model boundaries and whether human response is externally changed or affected by factors within the model. As a result, for models with endogenous representation of risk response, we identified two distinct categories based on how information diffusion is represented. Thus, three different categorizations of risk response are used: exogenous, endogenous (implicit information diffusion), and endogenous (explicit information diffusion). 
 
The explicit information diffusion models can be further separated into two formulations. The split-susceptible formulation is used in several papers \shortcite{abbas2022evolution,browne2022differential,costello2023model,epstein2021triple,kassa2020analysis}. The coupled information/behavior-disease formulation is used in a variety of other works \shortcite{juher2023saddle,li2022complex,Martcheva,morsky2023impact}. The model presented by \shortciteA{zhou2020effects} blends both types. \shortciteA{kumar2023multiscale} is an example of the second formulation as they present a neural model coupled with an SIR model. The neural model receives information about the disease in two forms: through disease information such as infections, deaths, and hospitalizations and through awareness campaigns. This stimulates firing of neurons in the brain in different regions depending on the type of information received. These combine to create a behavioral response, which is then incorporated into the epidemic model in the disease transmission term. Other models that incorporate this second formulation of explicit information diffusion include \citeA{juher2023saddle} (focus on awareness decay), \citeA{morsky2023impact} (focus on NPI adherence), and \citeA{Martcheva} (focus on compliance with social distancing policies).

The implicit information diffusion models, in contrast, tend to follow the same model type, generally SIR, possibly with additional infectious compartments. Disease information directly influences transmission without filtration through another source such as media or government and comes in the form of hospitalizations (e.g., \citeA{montefusco2022interacting}), deaths (e.g., \citeA{rahmandad2021behavioral,rahmandad2022enhancing,rahmandad2022quantifying}), and infectious population size (e.g., \citeA{alharbi2022nature,manrubia2022individual,rahmandad2021behavioral}). The functional forms of the implicit information diffusion risk response mechanism vary but often have a fractional form \shortcite{montefusco2022interacting,n2022effect,rahmandad2022enhancing}. The partial differential equation model from \citeA{song2022analysis} uses an exponential formulation of risk response affecting disease transmission with feedback coming from the infectious compartment but with a media function (although media is not a state variable) weighting the influence of prevalence. Furthermore, we note that not all models include differentiable functions for risk response. The model in \citeA{ihme2021modeling} is a good example as risk response can take two values depending on pandemic risks. Specifically, authors set a threshold in their model for per capita death (8 deaths per million population) after which infectivity rate drops, assuming that in such a scenario, social distancing has to be mandated, leading to change in contact rates.  

It is important to note that, while risk response is the most common behavioral mechanism modeled in behavioral epidemic models \cite{ferguson2007capturing, bauch2013social}, it is not the only phenomenon to be modeled. We found that several models layer these types of risk-response feedbacks with other human behaviors such as adherence fatigue (the tendency to avoid paying attention to risks over time) as in \shortciteA{de2021effect} or vaccination uptake (the tendency to vaccinate as a result of risk perception) as in \shortciteA{wang2023evaluating}, but these additional behavior considerations fall beyond the scope of this paper. 

The most common model structure observed in this compilation of compartmental models uses ordinary differential equations accompanied by linear or non-linear algebraic representations. Other types of model formulations use difference equations (e.g., \citeA{costello2023model}), stochastic systems (e.g., \citeA{manrubia2022individual, ochab2023multiple}), and partial differential equations (e.g., \citeA{song2022analysis}). Some authors also incorporate the capacity for different behavior by dividing the population based on characteristics such as age (e.g., \citeA{wirtz2021changing,dick2021covid,childs2022modeling}). 

\subsection{Information inputs to behavior change}

\begin{sidewaystable}
\small
\centering
\caption{\textbf{Information influencing transmission externally (exogenous models) or driving risk response internally (endogenous models).} Articles that appear in more than one box have more than one model with different characteristics.}
\begin{tabular}{|m{1in}|l|l|l|}
\hline 
& \textbf{Exogenous}  
& \begin{tabular}{@{}l@{}}\textbf{Endogenous} \\ \textbf{(Implicit Information Diffusion)}\end{tabular}
& \begin{tabular}{@{}l@{}}\textbf{Endogenous} \\ \textbf{(Explicit Information Diffusion)}\end{tabular}\\
\hline
\textbf{Disease \newline information:} \newline  
\textit{Infectious population size, number of deaths, hospital admissions} & \begin{tabular}{@{}l@{}}
\cite{alharbi2022nature}, Case 1
\end{tabular}
& 
\begin{tabular}{@{}l@{}}
\cite{alharbi2022nature}, Cases 2 \& 3 \\
\cite{ghaffarzadegan2021simulation} \\
\cite{ihme2021modeling} \\
\cite{lim2023similar}\\
\cite{manrubia2022individual}  \\
\cite{montefusco2022interacting} \\
\cite{n2022effect} \\
\cite{ochab2023multiple} \\
\cite{rahmandad2021behavioral} \\
\cite{rahmandad2022enhancing}\\
\cite{rahmandad2022quantifying} \\
\cite{song2022analysis} \\ 
\shortcite{zhang2023renewal} 
\end{tabular}
& 
\begin{tabular}{@{}l@{}} \cite{abbas2022evolution} \\
\cite{browne2022differential} \\
\cite{costello2023model} \\
\cite{juher2023saddle} \\
\cite{kassa2020analysis} \\
\cite{kumar2023multiscale} \\
\shortcite{macdonald2021modelling} \\
\cite{Martcheva} \\
\cite{morsky2023impact} \\
\cite{n2022effect} \\ 
\cite{wirtz2021changing}  
\end{tabular}
\\
\hline
\textbf{Disease \newline awareness (fear)} & & & 
\begin{tabular}{@{}l@{}}
\cite{abbas2022evolution} \\
\cite{browne2022differential} \\
\cite{epstein2021triple}  \\
\cite{kumar2023multiscale} 
\end{tabular}
\\
\hline
\textbf{Information index} & &  & 
\begin{tabular}{@{}l@{}}
\shortcite{bliman2022tiered} \\
\cite{buonomo2020effects} \\
\cite{d2022behavioral} \\
\cite{ochab2023multiple}  
\end{tabular}
\\
\hline
\textbf{Media \newline coverage} & & 
\cite{song2022analysis} & 
\begin{tabular}{@{}l@{}}
\cite{jasim2021study} \\
\cite{li2022complex} \\
\shortcite{rai2022impact} \\
\cite{zhou2020effects} 
\end{tabular}
\\
\hline
\textbf{Policy changes} & 
\begin{tabular}{@{}l@{}}
\cite{childs2022modeling}\\
\cite{dick2021covid}\\
\shortcite{eryarsoy2023models}
\end{tabular}& &  \\
\hline
\textbf{Age of \newline susceptible} & 
\begin{tabular}{@{}l@{}}
\cite{bubar2021model}\\
\cite{duggan2024age}
\end{tabular}& & \\
\hline
\textbf{Major events} & 
\begin{tabular}{@{}l@{}}
\shortcite{santos2021infusing}
\end{tabular}& & \\
\hline
\end{tabular}
\end{sidewaystable}
\label{table:diffusion}

In modeling risk response, a key decision involves choosing the information source that drives individuals to change their behavior. For example, do people react to a change in the number of cases, hospitalizations, deaths, or other factors about the pandemic? 
A second key decision is how to model the spread of information about the disease. Generally, information diffusion models are used to understand and predict how information spreads through a population or network. These models are crucial in various fields, including epidemiology, sociology, marketing, and technology management \cite{bass1969new}. Depending on the problem at hand and the importance of representing details about the spread of information, modelers may decide to more or less explicitly capture various mechanisms on the spread of information, such as word of mouth, media effects, or peer pressure. These models can be incorporated as a part of the epidemic model.

Table \ref{table:diffusion} shows models classified based on the type of risk response formulation and what information influences perceived risk and thus risk response decisions. Depending on the source of the information -- directly from the state of the disease or filtered through another source -- we categorize these models as implicit or explicit, respectively. The information can be directly related to the disease, such as the level of infection (prevalence), the number of new cases (incidence), the size of the quarantine class, or the number of deaths. In such scenarios, change in behavior is directly influenced by disease information, usually via implicit information diffusion. When disease information reaches the public directly, some models (e.g., \citeA{rahmandad2022enhancing}) show another type of filtration that occurs in the form of a time delay (see Section \ref{delays}). We still consider this formulation to be implicit information diffusion as the behavior change is tied directly to the disease state, just with a delay such that disease information does not reach individuals instantaneously.

In contrast, as stated the information can also be filtered through another source such as media coverage, social media, government news dispersion, or an information index (see below for definition) before it reaches the general population. Categorized as explicit information diffusion, the information filter is another state variable and is dependent upon disease information (e.g., level of infection, number of new cases, number of deaths, etc.). The change in behavior is then dependent upon the information available, after passing through the filter considered such as media coverage or the government. \citeA{bliman2022tiered}, \citeA{buonomo2020effects}, \citeA{d2022behavioral}, and \citeA{ochab2023multiple} 
all consider the information index as the driving risk response mechanism. An \textit{information index} is a risk response function summarizing the disease information reaching the general public that aids in determination of the perceived risk. Typically, an information index function includes a distributed delay, as memories of past disease levels remain relevant to inform behavior \cite{buonomo2020effects}.

\section{Examples of risk response formulations}

To give a clearer picture of categories of model formulation for risk response, we highlight particular models. We focus on the key components related to the inclusion of risk response and  discuss how the behavior mechanism of risk response impacts conclusions from the studies. 

\subsection{Classic compartmental epidemic modeling}

The modeling of the spread of infectious disease through a population has a long history, dating back more than a hundred years \cite{kermack1927contribution,ross1916application,ross1917applicationII,ross1917applicationIII} and is classically done by dividing the population based on disease characteristics. The most general being the SIR model where individuals fall into a single type -- susceptible (S), infectious (I) or removed/recovered (R) -- and move between compartments as the result of infection (S to I) or recovery (I to R). Importantly, the infection term involves the interaction between susceptible individuals and infectious individuals, which is inherently captured by a non-linear mathematical formulation. The elegance of these models lies in their simplicity as well as their relevancy in producing characteristic epidemic curves. Furthermore, they are easily extendable through the addition of compartments representing other disease states, such as latently infected. In that case, individuals that are infected but not yet infectious are separated in what is standardly represented by E for exposed producing the SEIR model. Both the SIR model and its most basic extension, the SEIR model, are widely used, studied, and taught (e.g.  \citeA{anderson1991infectious,anderson1979population,keeling2011modeling,kermack1927contribution,vynnycky2010introduction}).

\subsection{Exogenous risk response}

The classic SIR model (and its extension, the SEIR model) assumes a single constant transmission rate parameter which combines susceptibility, infectivity, and contact structure, all of which may be changing through the course of an epidemic \cite{kermack1927contribution}. One simple way to account for these behavioral adaptations is to use different constant transmission rates for fixed intervals.  Often the length and timing of these intervals is determined by some alternative information source, such as governmental policy. When fitting a model to data, the history of the external drivers, such as the level of government policy, is known and thus can be incorporated. 

As an example, consider the model by \citeA{dick2021covid} on waning immunity and vaccination in Canada. The authors fitted the transmission rate over windows chosen to coincide with changes to governmental policy in Canada. As contact structure in their age-structured model also changed at the same time points, transmission rate encompassed changes to viral infectivity and compliance to NPIs. While they could tune their model to past data by the frequency and length of the windows, for all forecasting they had to rely on predefined future scenarios about contact structure and government policy. In summary, their exogenous structure could estimate multiple past waves of the pandemic due to exogenously changing the level of transmission as the country changed policy, but for accurate future projections they would need data on future policies and public compliance with those policies. This is because the exogenous formulation does not allow for feedback between the behavior and disease dynamics.

Another example of exogenous modeling of risk response is by \citeA{bubar2021model}. In their compartmental model of vaccine prioritization, average contact rate is constant over time but different for different sub-populations. Analyzing a variety of vaccination policies, they assume that people do not increase their contact rate as the number of cases or death declines. This can potentially lead to over-estimation of the impact of vaccination.   

We note several examples of exogenous formulation of infectivity in the literature that deal with optimization, common in Operations Research \cite{currie2020simulation}. In this body of literature, assuming public compliance with policies is constant, compartmental models are formulated with time-varying inactivity and modelers seek values of infectivity that minimizes total costs \cite{eryarsoy2023models}. The concern, however, is that in this approach public risk response is limited to what is enforced by the government.

\begin{figure}[ht!]
\centering
\includegraphics[width=.4\textwidth]{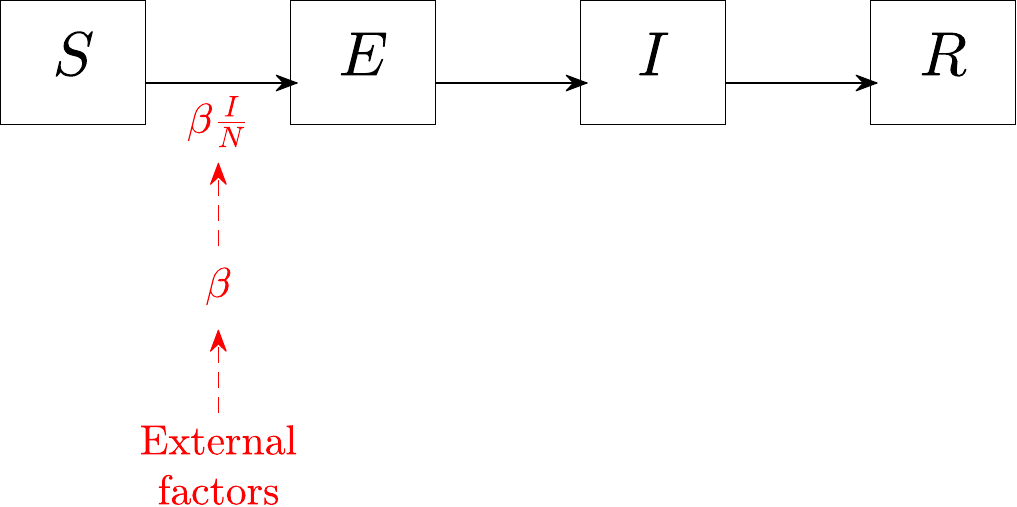}
\caption{Flow diagram describing the general structure of exogenous SEIR models. Red dashed lines indicate influences on disease transmission. Solid lines denote movement of individuals. The compartments are: $S$, susceptible; $E$, exposed; $I$, infectious; $R$, removed.}
\label{fig:exogenous}
\end{figure}

Overall, models with an exogenous representation of risk response can better replicate past data and provide precise estimations of changes in contact rates or population flows between different compartments if proper parameter estimation is performed \cite{duggan2024age}. A constant parameter for infectivity or contact rate in the conventional SEIR format will only generate one wave of the disease on an epidemic timescale. The alternative, used within exogenous modeling, is to let the related parameter(s) change through time, for example by setting different values at different time intervals. This approach can help offer a more consistent narrative of past trends especially when specific policies were implemented or the public changed their behavior. However, when using contact rate or infectivity as external inputs to the model, model-based explanations for why these parameters are changing are lost (see flow diagram in Figure \ref{fig:exogenous}). Additionally, since future data on these inputs are unavailable (e.g., we cannot know what will be future contact rates), the accuracy of forecasts and policy analyses would depend on selecting a proper scenario of future contacts and infectivity. In most cases, however, modelers tend to assume that people's response to the pandemic will remain unchanged. 

\subsection{Endogenous risk response with implicit information diffusion}

\begin{figure}[ht!]
\centering
\includegraphics[width=.4\textwidth]{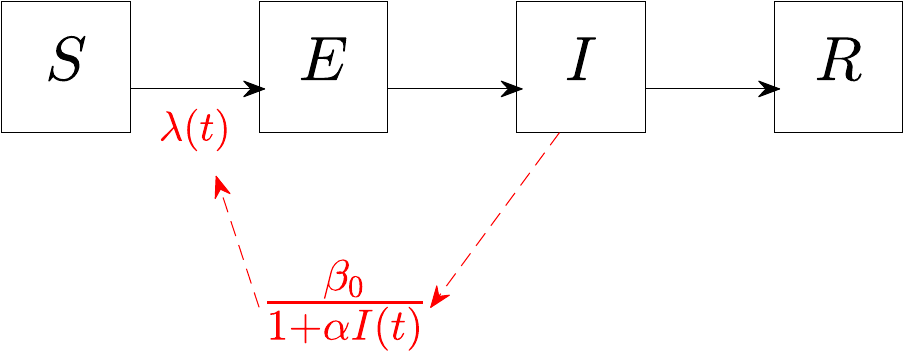}
\caption{Flow diagram describing the model from N'konzi et al. (2022).  Influence of behavior changes occurs endogenously via implicit information diffusion. Red color indicates contributions to risk response formulation. Solid lines denote movement of individuals, while dashed lines indicate information feedback. Compartments are  $S$, susceptible; $E$, exposed; $I$, infectious; $R$, recovered. }
\label{fig:nkonzi}
\end{figure}

In this approach, not only is the assumption about time-independent infectivity relaxed, but infectivity is formulated as a function of the state of the disease. This closes a balancing (negative) feedback loop in which infectivity increases transmission and number of cases, which with a time lag, decreases infectivity and transmission rate. \citeA{n2022effect} construct a model of this type to analyze the spread of COVID-19, where the force of infection is given by $\lambda(t)=\frac{\beta_0 I(t)}{1+\alpha I(t)}$, with $I(t)$ the ratio of the number of infected people to the total population size. The transmission rate is given by the fraction $\frac{\beta_0}{1+\alpha I(t)}$ such that as $I(t)$ increases, the transmission rate declines. Here, $\beta_0$ represents the effective contact rate, without disease control measures, which may depend on demographics and density of a region, and $\alpha$ represents the level of government policy response to control the epidemic. 
See Figure \ref{fig:nkonzi} for the flow diagram of individuals and note the feedback loop.

Several studies considered in this review analyze models of a similar structure; however, the behavioral feedback loop (the dashed line in the figure) can come from daily infection, positive rates of tests, or daily death. For example, the paper by \citeA{rahmandad2022enhancing} formulates $\beta$ as an inverse function of lagged death rate. Due to their relatively simple structure, these models have only a few extra parameters in comparison to a conventional SEIR model, such as $\alpha$ or time delay between the state of the disease and its effect on infectivity, which can be estimated through model calibration. Many of these SEIRb models were calibrated with various data sources including death rate in different US states for the purpose of forecasting \cite{rahmandad2022enhancing}. Another study uses data from more than 100 countries, fitting a similar model, for the purpose of examining responsiveness. This parameter turns out to significantly explain the variation in death rate across countries \cite{lim2023similar}.  

\subsection{Endogenous risk response with explicit information diffusion}

The explicit information diffusion models can be further separated into two formulations: the split-susceptible formulation and the coupled information/behavior-disease formulation.

\subsubsection{Split-susceptible formulation}

\begin{figure}[h!]
\centering
\includegraphics[width=.5\textwidth]{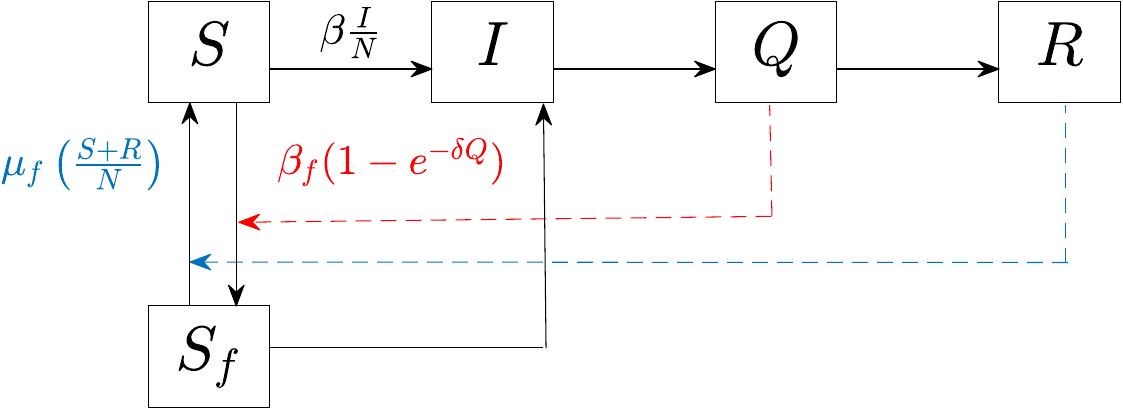}
\caption{Flow diagram describing the model from Abbas et al. 2022).  Red and blue coloring indicate contributions to risk response formulation with blue coloring for fear relaxation. Solid lines denote movement of individuals, while dashed lines indicate influence. Influence of behavior changes occurs endogenously via explicit information diffusion with split susceptible compartment. Compartments are $S$, susceptible; $S_f$, fearful susceptible; $I$, infectious; $Q$, quarantined; $R$, recovered.}
\label{fig:abbas}
\end{figure}

Some endogenous models explicitly separate sub-populations based on their risk perceptions. \citeA{abbas2022evolution} is an example of the first type of explicit information diffusion model where the susceptible population is split. In particular, the susceptible population is divided based on the level of fear between standard susceptibles, $S$, and those that are aware or fearful, $S_f$.  Risk response is included as responding to fear of the disease by adhering to NPIs, i.e., social distancing, avoiding social gatherings, and mask wearing. 
Fearful susceptible individuals in $S_f$ change their behavior in response to increased perceived risk, and therefore experience a reduced rate of disease transmission compared to individuals in $S$. The movement into $S_f$ is determined by the population that is currently quarantined, $Q$, which leads to the change in the transition from $S$ to $S_f$, given by $\beta_fS(1-e^{-\delta Q})$ where $\beta_f$ is the rate of transmission of fear and $1/\delta$ is the average number of confirmed cases reported by news. The return of fearful individuals to full susceptibility depends on the fraction of the population that is recovered, $R$, or fully susceptible, $S$, given by $\mu_f S_f \left(\frac{S+R}{N}\right)$, where $\mu_f$ is the rate of relaxation of behavioral change. See Figure \ref{fig:abbas} for flow diagram of individuals.

\citeA{abbas2022evolution} examine case studies in South Korea, Pakistan, and Japan, where the model is fit to data for each of the three countries over a time period which included three outbreak waves. Simulations over the period of about a year show the active cases decreasing in each country; with lower transmission rate $\beta$, the outbreak peak occurs later and is smaller. Only for Japan do the simulations show small waves occurring after the initial outbreak in the simulations. Additional analyses show that the rate of behavioral change relaxation, $\mu_f$, influences outbreak peak timing and intensity, and $\beta_f$, the rate of transmission of fear has a stronger influence on peak size. In each successive wave, the interval of sensitivity to the combination of $\beta_f$ and $\mu_f$ becomes larger; in other words, the populations require more diligence in NPI usage and more transmission of fear for each successive wave to lower the peak in cases, likely due to public sensitivity to the disease waning over time. Thus, maintaining behavior change through $\beta_f$ and $\mu_f$ is critical to reduce successive peak sizes, requiring adjustments in policy to maintain low caseloads for each successive outbreak.

\subsubsection{Coupled information/behavior-disease model}
\label{sec:endo:exp:coup}

\begin{figure}[ht!]
\centering
\includegraphics[width=.4\textwidth]{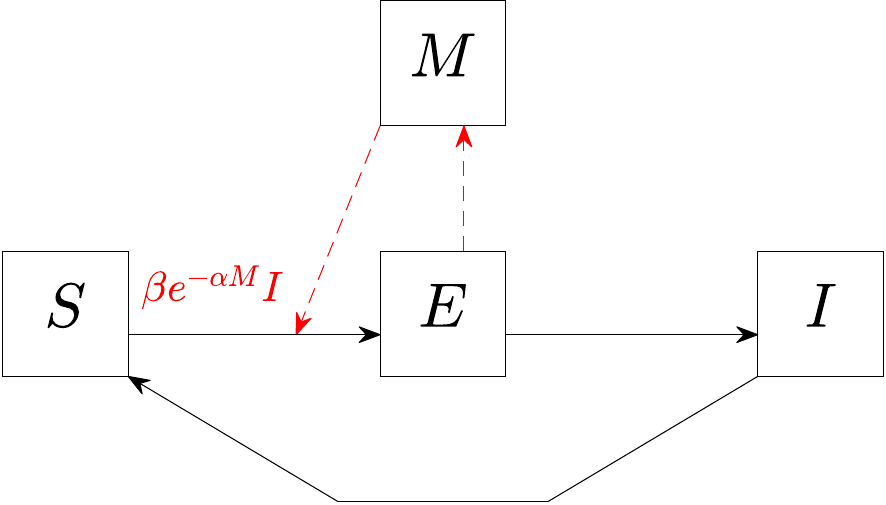}
\caption{Flow diagram describing the model from Li \& Xiao (2022). Red coloring indicates contributions to risk response formulation. Solid lines denote movement of individuals, while dashed lines indicate information feedback. Influence of behavior changes occurs endogenously via explicit information diffusion. Compartments are $S$, susceptible; $E$, exposed; $I$, infectious; $M$, media.}
\label{fig:li}
\end{figure}

\citeA{li2022complex} is an example of the second type of explicit information diffusion model with information spread via explicit information diffusion. Information about the disease diffuses through an additional state variable $M$ which represents media coverage of the disease and is coupled with an SEI disease component. The disease transmission rate is dependent upon the function $e^{-\alpha M}$, where $\alpha$ is the weight of the media on the contact rate. In turn, growth in media coverage $M$ positively depends upon $E$ (exposed individuals), closing the feedback loop in the coupled model. See Figure \ref{fig:li} for the flow diagram of individuals. In contrast to other models highlighted, \citeA{li2022complex} do not consider immunity and thus do not include a removed compartment; rather, individuals who recover move directly back to the susceptible compartment. Another way of taking the same approach is to represent information diffusion with a state variable of the fraction of aware individuals (or fraction adopting NPIs) and link this back to the transmission rate (e.g. \citeA{Martcheva}).  

Numerical simulations show that the model is sensitive to, in decreasing order, disease transmission rate ($\beta$), rate of progression from exposed to infectious, recovery rate, weight of the media on the contact rate ($\alpha$), reporting rate, saturation effects of media reports, and the effect of medical resource limitation,  in terms of the total number of infections after one month of an outbreak. Additionally, simulations show the existence of multiple epidemic waves for the infectious population, and analysis suggests that oscillatory dynamics of the model are modulated by the nonlinear terms driving recovery and media coverage.

\section{Mathematical Formulations}
\label{sec:math}

In this section, we construct models based on the incorporation of discussed mechanisms from the literature. We perform simulations and evaluate how the different risk response mechanisms affect outbreak size and timing, as well as the potential for multiple waves.

\subsection{Base model}

\begin{figure}
\centering
    \includegraphics[width=.5\textwidth]{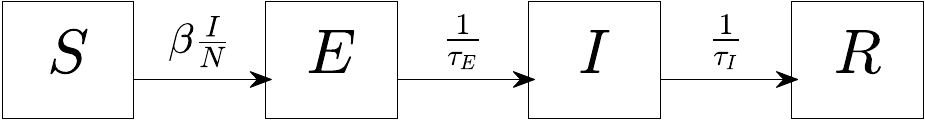} 
\begin{tabular}{cc}
 \begin{minipage}{0.5\textwidth}
 \vspace{0.5cm}
 \begin{tabular}{cl}
\hline Variables & \\
\hline $S$ & number of susceptible individuals \\
$E$ & number of infected (not infectious) individuals \\
$I$ & number of infectious individuals \\
$R$ & number of removed individuals \\
$N$ & total population size $(=S+E+I+R)$ \\
\hline Parameters & \\
\hline$\beta(\cdot)$ & transmission rate \\
$\tau_E$ & time from infection to infectiousness \\
$\tau_I$ & time from infectiousness onset to recovery \\
\hline
\end{tabular}
 \end{minipage}&
 \begin{minipage}{0.5\textwidth}
 $$
 \begin{aligned}
 & \frac{d S}{d t}=-\beta(\cdot) S \frac{I}{N} \\
 & \frac{d E}{d t}=\beta(\cdot) S \frac{I}{N}-\frac{E}{\tau_E} \\
 & \frac{d I}{d t}=\frac{E}{\tau_E}-\frac{I}{\tau_I} \\
 & \frac{d R}{d t}=\frac{I}{\tau_I}
 \end{aligned}
 $$
 \end{minipage}
\end{tabular}
    \caption{Classical SEIR model without demographics. \textit{Top:} Flow diagram of individuals through compartments $S$, susceptible; $E$, exposed/infected but not infectious; $I$, infectious; $R$, removed. \textit{Bottom Left:} Variables and parameters. \textit{Bottom Right:} System of differential equations for the base model. }
    \label{fig:SEIRmodel}
\end{figure}

A simplified epidemic model is described by the system given in Figure \ref{fig:SEIRmodel}. The population is divided by disease status into four compartments: susceptible ($S$), infected but not infectious ($E$), infectious ($I$), and removed ($R$). As we assume there is no movement of individuals into or out of the population via death, birth, or migration, the total population size is constant with $N = S+E+I+R$. The susceptible population is infected at frequency dependent rate $\beta(\cdot) \frac{I}{N}$, where $\beta(\cdot)$ is the infectivity rate. We denote infectivity rate by $\beta(\cdot)$ as it is constant in the most basic version but may be a function of other states of the system in more complex formulations. The infected population becomes infectious at rate $1/\tau_E$, and the infectious population recovers at rate $1/\tau_I$. Thus, the average duration in state $E$ is $\tau_E$ and in state $I$ is $\tau_I$. 

\subsection{Exogenous formulation of infectivity}

To highlight the difference between exogenous and endogenous formulations, we first consider the SEIR model with an exogenous behavior incorporation, a standard formulation in epidemic modeling \cite{kermack1927contribution,anderson1979population,martcheva2015introduction}. The most basic exogenous formulation uses constant infectivity, i.e., $\beta(\cdot)=\beta_0$. Note that here, we have added the subscript to distinguish it from the function $\beta$, although many models refer to a constant $\beta$. 

When taking into account seasonality or travel trends such as increased mobility during holidays, transmission rates in a population change over time. Thus, transmission can be written as a function of time, for example, $\beta(t)=\beta_0(1-kt)$, where transmission is a decreasing function of time until it reaches zero. This decrease in transmission rate could occur as policies are implemented to limit spread in a community.

\begin{figure}[hb!]
    \centering
\begin{subfigure}[]{.33\textwidth}
\centering
    \includegraphics[width=\textwidth]{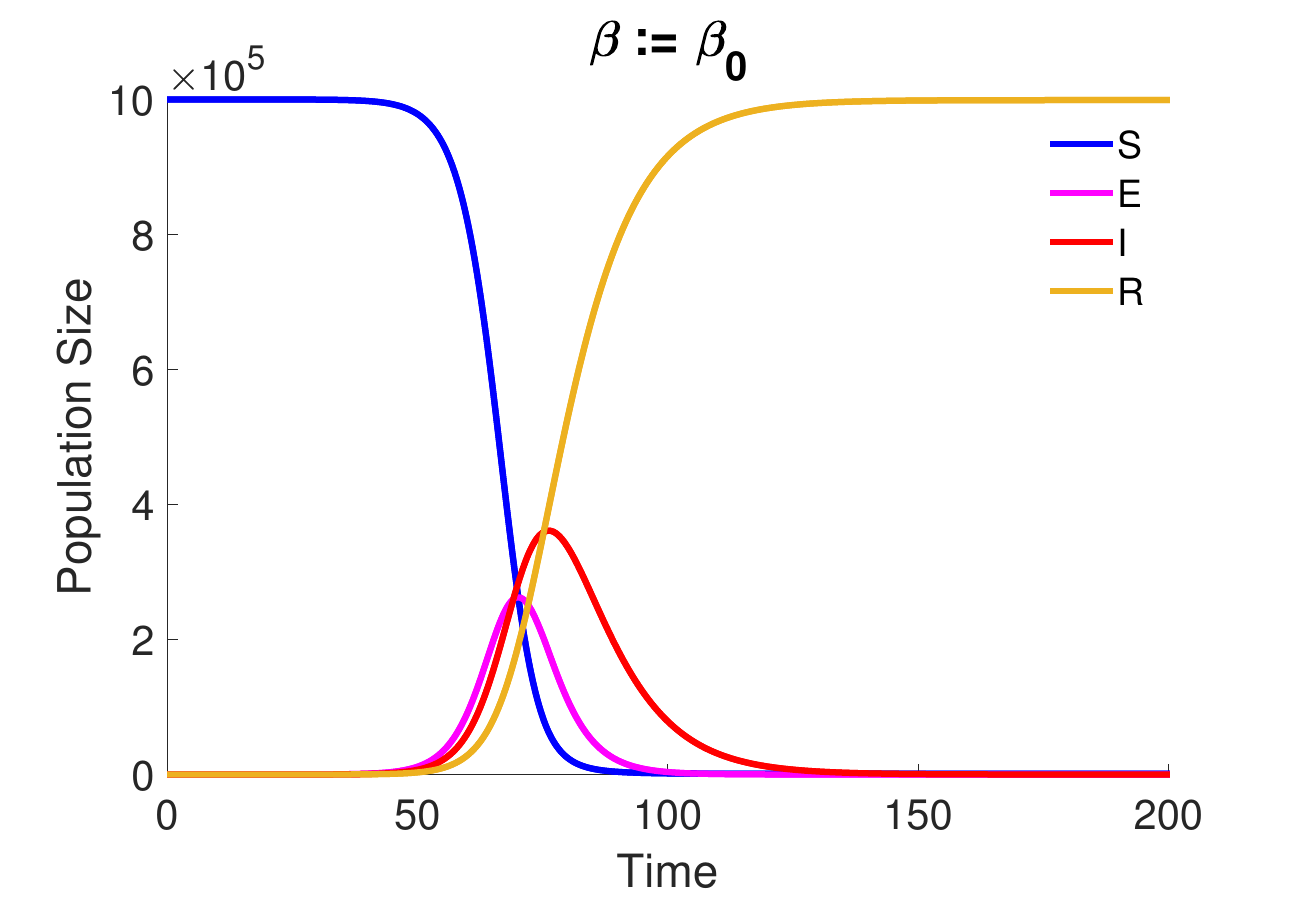}
    \caption{Constant transmission}
    \label{fig:exog:constant}
\end{subfigure}%
~
\begin{subfigure}[]{.33\textwidth}
\centering
    \includegraphics[width=\textwidth]{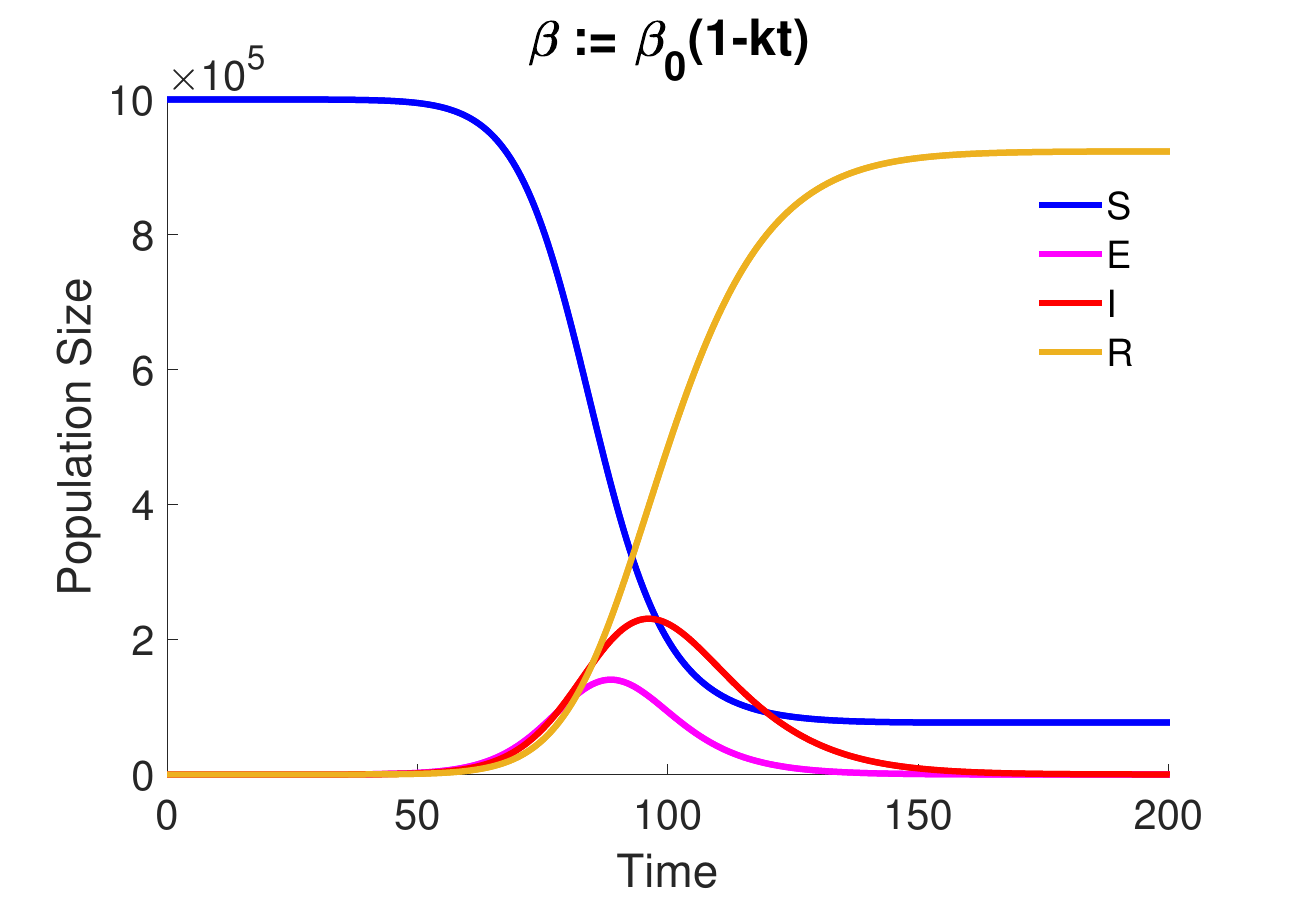}
    \caption{Time-varying transmission}
    \label{fig:exog:timevary}
\end{subfigure}%
~
\begin{subfigure}[]{.33\textwidth}
\centering
    \includegraphics[width=\textwidth]{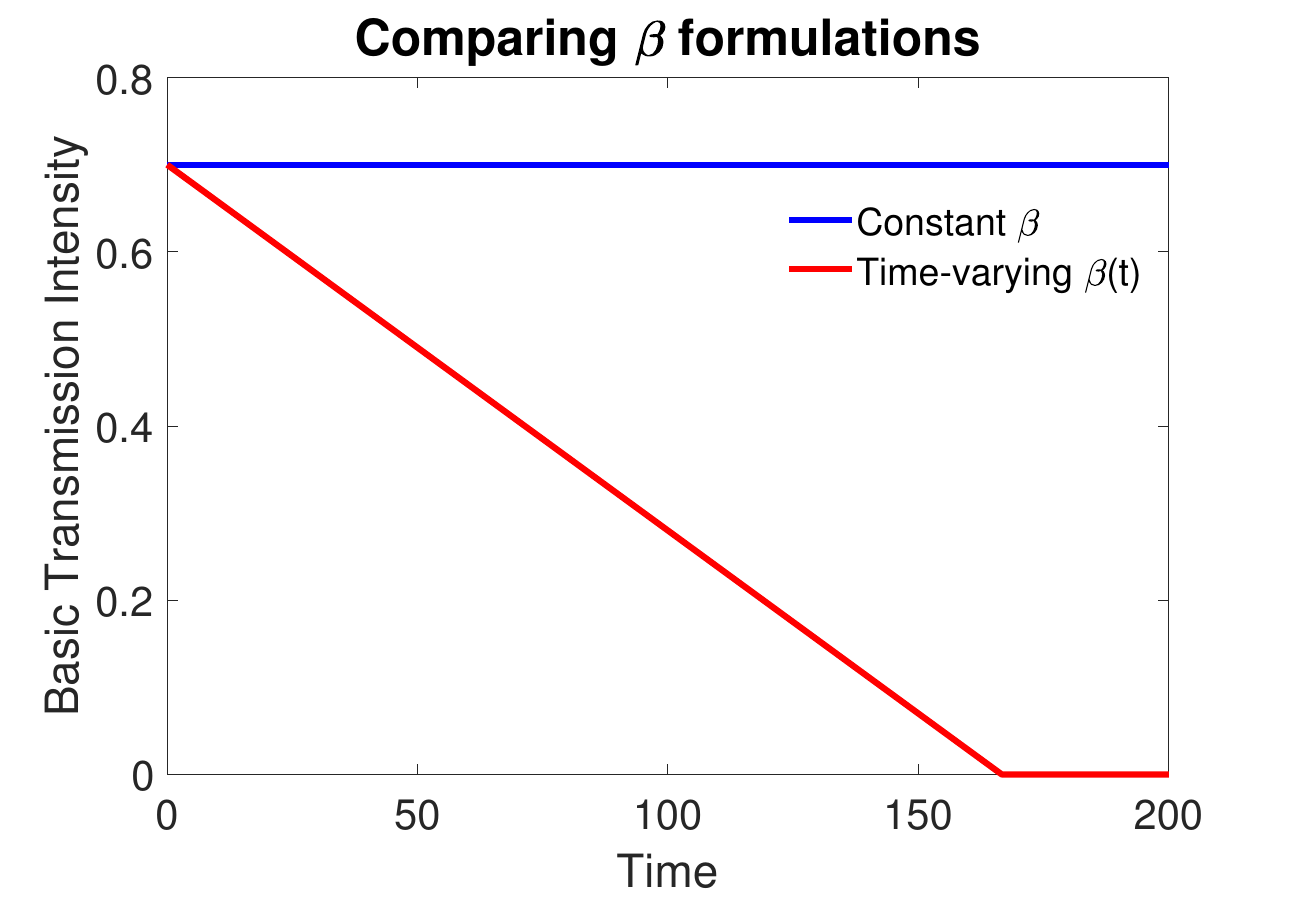}
    \caption{Transmission function}
    \label{fig:exog:trans}
\end{subfigure}
\caption{Comparison of exogenous infectivity formulations. Epidemics with (a) constant transmission and (b) time-varying transmission with transmission decreasing through time. The values of transmission rate through time are given in (c). Note that time-varying transmission is zero when $t>\frac{1}{k}$. Parameter values are listed in Table \ref{table:params} in the Appendix.}
\label{fig:exog}
\end{figure}

We compare outbreaks from constant transmission and time-varying transmission in Figure \ref{fig:exog}. Both transmission settings show similar dynamics: a single outbreak occurs before the system stabilizes without an infectious population. The outbreak in the chosen variable transmission model has a smaller magnitude (as expected, since the transmission rate is decreasing from the constant transmission setting as seen in Figure \ref{fig:exog:trans}), and the susceptible population has a larger steady state compared to the constant transmission model. As both show just a single epidemic wave before the disease decays to zero, these both are relatively unrealistic models for many real world infectious disease outbreaks, such as COVID-19. While the time-varying transmission model takes into account fixed changes in transmission over time, neither model takes into consideration how changes in disease severity affects transmission. 

\subsection{Endogenous formulation of infectivity}

We consider three different endogenous formulations to include behavior in the model as a modifier of infectivity. In these examples, the infectivity parameter $\beta(\cdot)$ is a function of some information cues about the state of the disease. We intentionally keep the models simple to focus on different formulations of behavioral risk response and compare such formulations. For example, the models we describe do not consider other important factors such as asymptomatic infectious cases, seasonality, hospitalization, or waning immunity.

\begin{table}[h]
\caption{Mathematical formulations of risk response. Formulations given here include a first-order time delay. In the case of no delay, the $F$ in $\beta(\cdot)$ is replaced by $I$.}
    \centering
    
    \begin{tabular}{lcll}
    \hline
         Type &     $\beta(\cdot)$ &  Additional variables and parameters \\
         \hline
         \hline
         \multicolumn{2}{c}{\textbf{Exogenous}}\\
         \hline
         
         Constant & $\beta_0$ \\
         \hline
         Time-varying & $\left\{
     \begin{array}{lr}
     \beta_0(1-kt), & 0 \leq t \leq \frac{1}{k},\\
     0, & t > \frac{1}{k},
     \end{array}\right. $
     & $k$, time-dependent decay\\
         \hline
         \hline
         \multicolumn{2}{c}{\textbf{Endogenous}}\\
                   \hline
         \multirow{3}{*}{ Time-delay equation} & \multirow{3}{*}{$\displaystyle\frac{dF}{dt} = \frac{I-F}{\tau_F}$} & $F$, perceived number of caseloads \\
         & & $\tau_F$, delay between infectiousness \\
         & & $\qquad$ and information dissemination\\
         \hline
          \multirow{2}{*}{(1) Fractional} &  \multirow{2}{*}{$\displaystyle\frac{\beta_0}{(1+\alpha F)^{\gamma}}$} 
          & $\alpha$, community sensitivity to prevalence \\
          & & $\gamma$, diminishing impact of infection risk\\
          \hline
          (2) Exponential & $\beta_0 \exp(-kF)$ & $k$, proportionality constant \\[.2cm]
          \hline
          \multirow{3}{*}{(3) Explicit} & \multirow{3}{*}{$\beta_0 (1-X)$, where}  & $X$, proportion of individuals adhering \\
          & & $m$, infectivity multiplier \\
          & $\displaystyle\frac{dX}{dt} = X(1-X)(mF-c)$ & $c$, cost of preventative measures \\
          & & & \\
         \hline
          \hline
    \end{tabular}
    \label{tab:formulations}
\end{table}

\subsubsection{Implicit information diffusion}

We focus on two examples of implicit information diffusion: a fractional formulation and an exponential formulation. In both, transmission is structured as a decreasing function of infection prevalence. As infections increase, transmission  decreases, capturing human behavior, such as an increase in protective measures leading to a decrease in risky behavior. As transmission decreases, infection levels decrease, which leads to a relaxation of protective measures and thus an increase in transmission and an increase in infections. This feedback loop results in oscillating infection levels. 

Basic formulations of each type (which each include a first-order time delay) are found in Table \ref{tab:formulations}. In the fractional formation, as found in \citeA{rahmandad2022enhancing}, infectivity becomes inversely related to the number of perceived cases.  
The exponential formulation, as seen in \citeA{ghaffarzadegan2021simulation}, considers an exponential structure for infectivity related to the number of perceived infections, a first-order time delay for the occurrence of infections. In both instances, risk response is driven by perceived risk of infection, measured by perceived infections with a time delay between infections occurring and public awareness of infection level.

\subsubsection{Explicit information diffusion}

Explicit information diffusion occurs either in the form of a split-susceptible formulation or a coupled information/behavior-disease formulation. Individual susceptibility changes either by movement between fully susceptible and less susceptible compartments or by a change in the disease transmission rate. In either case, behavior change is driven by information such as in the form of incidence, media, information index, or government intervention. 

Our explicit information diffusion formulation is based on individuals choosing optimal strategies for behavior as in \citeA{Martcheva}. In addition to the disease classes, there is now a state equation for $X$, the proportion of individuals adhering to a strategy, for example, a particular NPI such as mask wearing. In \citeA{Martcheva}, $X$ represents social-distancing. The infectivity rate ($\beta_0$) is constant in the absence of behavior change, but the transmission term is affected by the proportion of individuals not adhering to the given intervention. Here, the time delay is incorporated in the behavior sub-model, rather than in the transmission term (as in the implicit information diffusion model). In the case shown in Table \ref{tab:formulations}, the proportion of individuals adhering is dependent on perceived caseloads.

\subsubsection{Incorporation of delays}\label{delays}

Realistically, individuals do not have instantaneous information about changing case levels, such that there is a delay between when these events occur, when they are reported, and when individuals receive and act on this information. 
Thus, risk response is typically constructed based on perceived risk from information about the state of the disease that is several hours (or days or weeks) old and does not accurately reflect the current state of the epidemic.

We consider different time delays and their impact on transmission by building layers of information lags into our model as adding or reducing the number of delays changes epidemic behavior. For example, one time delay could arise due to a lag between when infections are treated in a hospital and reporting of those infections to a government agency (first order delay). An additional delay could represent a lag between when the agency receives the data and when the data is released to the media (second order delay). A third delay could occur due to a lag between when the media receives the data and when the public is made aware of the information (third order delay).  

We build time delays into our model through the incorporation of additional equations. The time delay equation given in Table \ref{tab:formulations} is added on to each respective model for a first order time delay. For each subsequent $n$th order time delay, using the implicit information diffusion exponential model as an example,  additional equations are added of the form
\begin{equation*}
\frac{dF_i}{dt} = \frac{F_{i-1}-F_i}{\tau_F/n}
\end{equation*}
for $i \in\{2,\cdots,n\}$ and $\beta(\cdot)$ becomes a function of $F_n$. Additional time delays are added to the implicit information diffusion fractional model and explicit information diffusion models in a similar manner.

Other types of time delays are possible and are used in the models listed in Section \ref{sec:formulationsriskresponse}. Models using an information index (e.g., \citeA{ochab2023multiple,bliman2022tiered,buonomo2020effects,d2022behavioral}) all incorporate a time delay. \citeA{zhang2023renewal} formulates disease transmission using a fractional function, which includes a fixed delay, equivalent to an infinite order time delay of our formulation.

\subsubsection{Simulation results}

\begin{figure}[ht!]
 \centering
\begin{subfigure}[]{.3\textwidth}
\centering
    \includegraphics[width=\textwidth]{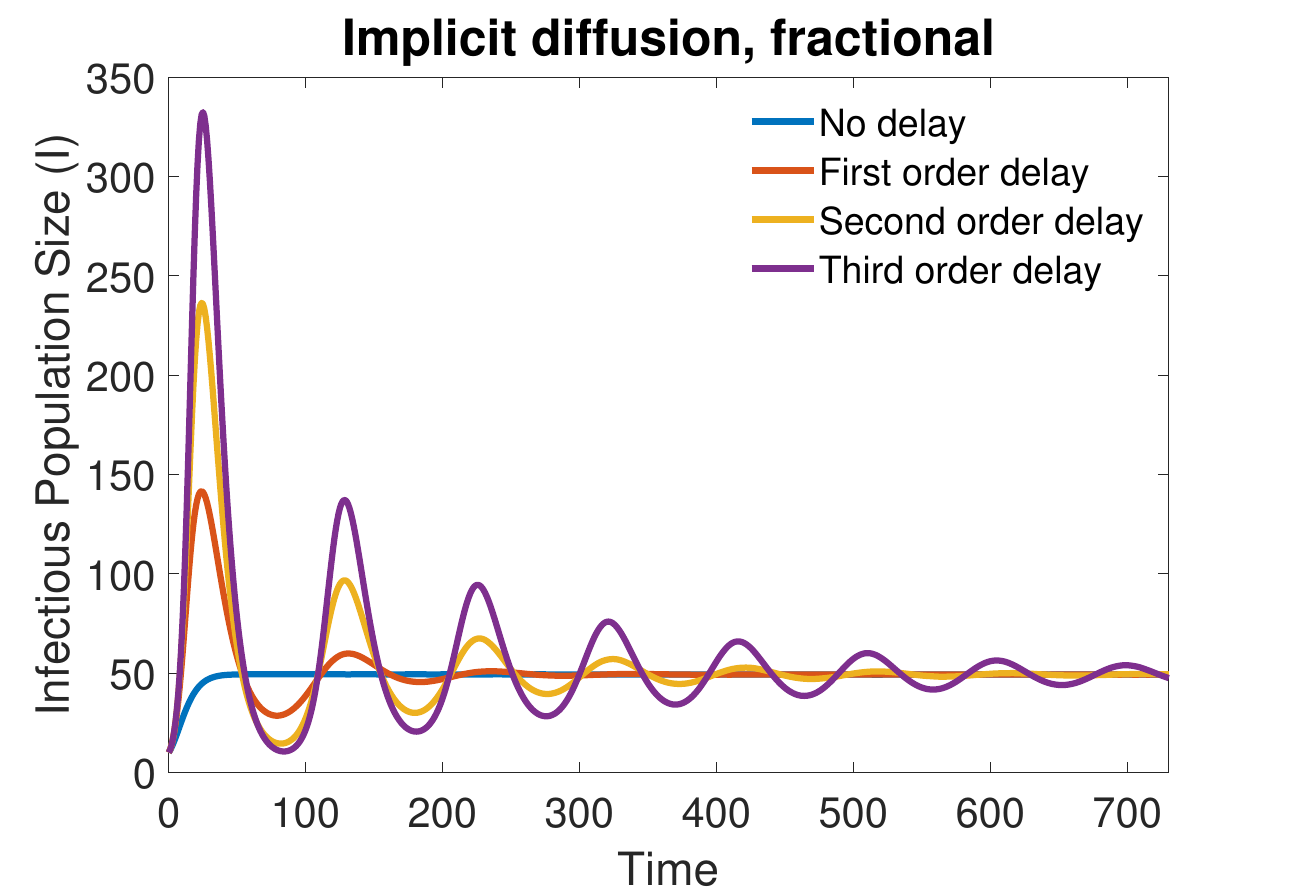}
    \caption{Implicit, fractional}
    \label{fig:endo:frac}
\end{subfigure}
\begin{subfigure}[]{.3\textwidth}
\centering
    \includegraphics[width=\textwidth]{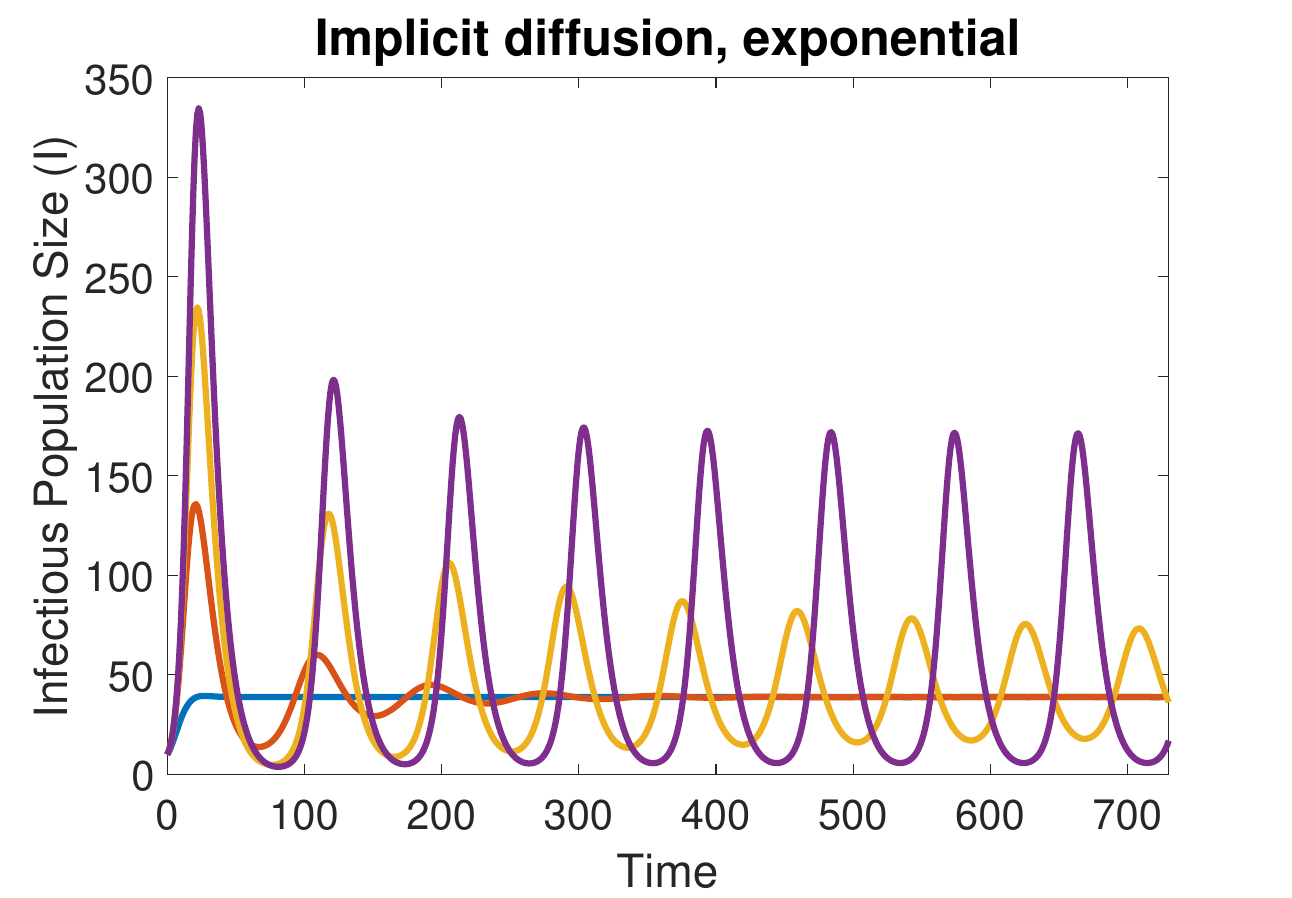}
    \caption{Implicit, exponential}
    \label{fig:endo:exp}
\end{subfigure}
\begin{subfigure}[]{.3\textwidth}
\centering
    \includegraphics[width=\textwidth]{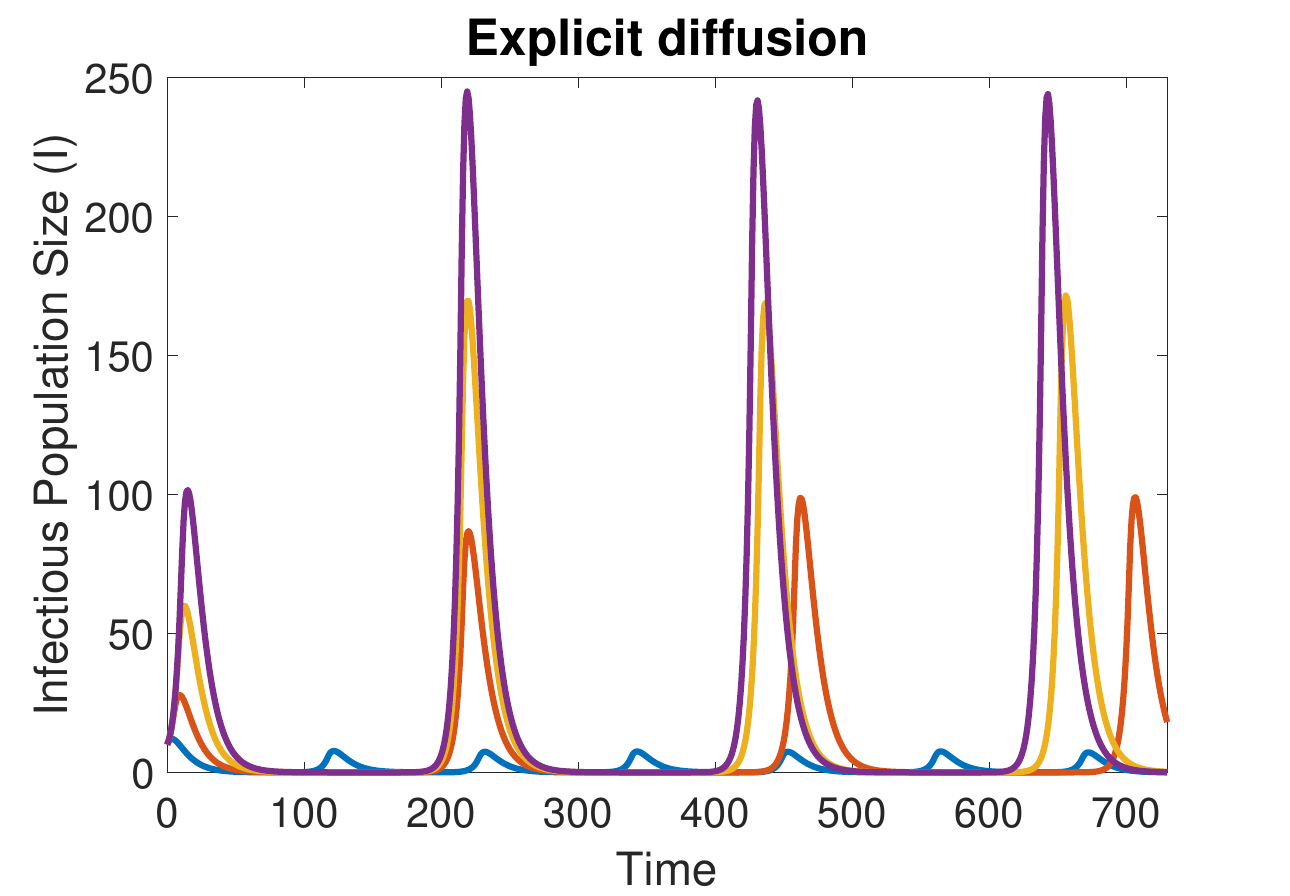}
    \caption{Explicit, SEIR\textit{x} model}
    \label{fig:endo:explicit}
\end{subfigure}
\caption{Comparison of endogenous formulations of risk response. Behavior is incorporated via (a) implicit information diffusion with fractional form, (b) implicit information diffusion with exponential form, and (c) explicit information diffusion. Line color in (a), (b) and (c) denote the order of delay. Parameters are given in Table \ref{table:params} in the Appendix.
}
\label{fig:endo}
\end{figure}

We examine the simulations of the three types of endogenous formulations with different orders of time delays (Figure \ref{fig:endo}). In each graph, the $x$-axis represents time given in days, and the $y$-axis represents the number of infectious individuals $I$ at that time. Figures \ref{fig:endo:frac} and \ref{fig:endo:exp} are simulations using the SEIR structure, as described in Figure \ref{fig:SEIRmodel}, modified by the respective fractional (1) or exponential (2) formulation for $\beta(\cdot)$ as given in Table \ref{tab:formulations}. These model structures assume permanent immunity to the disease upon recovery (i.e. there is no flow from $R$ back to $S$). Figure \ref{fig:endo:explicit} represents a structure for information diffusion coupled with the infectious disease component using the explicit (3) formulation for $\beta(\cdot)$ as given in Table \ref{tab:formulations}. Despite the flow of information, there is no flow of recovered individuals back to susceptible.

With implicit information diffusion where transmission $\beta$ is a fractional formulation (Figure \ref{fig:endo:frac}), increasing the number of delays increases the amplitude, but not the frequency, of oscillations. Additionally, the length of time before oscillations dampen to equilibrium increases with the number of delays. With implicit information diffusion where transmission $\beta$ is an exponential formulation (Figure \ref{fig:endo:exp}), the amplitude of the oscillations and the time until dampening to equilibrium increases, but the frequency of oscillations decreases. Both the fractional formulation and the exponential formulation with no time delay appear to reach an endemic equilibrium, although long-time simulations show that these will eventually reach zero over time (not shown). With implicit information diffusion with no time delay (Figures \ref{fig:endo:frac} and \ref{fig:endo:exp}), the population of infectious individuals quickly reaches equilibrium, and increasing the order of delay creates more outbreak waves. Exponential formulation (Figure \ref{fig:endo:exp}) exhibits the most consistent oscillatory dynamics, especially as the order of the delay increases.

With explicit information diffusion where infectivity $\beta_0$ is constant and the transmission term is affected by the proportion of individuals not adhering to NPIs given by state variable $X$, as in Figure \ref{fig:endo:explicit}, the oscillation frequency is less than in the implicit information diffusion models. Furthermore, the model requires high reactivity to infections from the population in order for the peak size to be a similar magnitude as in the implicit information diffusion models. Varying the cost of adherence ($c$) changes the frequency of oscillations: as $c$ decreases, oscillations appear less frequently. The simulation with no time delay continues to oscillate but with a small maximum peak size ($\leq 5$). All time delays show a lower frequency than the simulation with no time delay. However, as the order of the delay increases, oscillation frequency and maximum peak size increase.

All three formulations display outbreak waves over the course of a two-year epidemic, demonstrating that endogenous inclusion of the human behavior feedback loop allows epidemic models to exhibit more representative dynamics of disease outbreaks. As the order of time delays increases, peak size increases for all waves in the epidemic. With no time delay, only the explicit information diffusion model exhibits outbreak waves.

\section{Discussion}

The objective of this study was to examine various ways of modeling human behavior in epidemics with a particular focus on how people’s response to risks changes during disease outbreak. This issue has become especially important after the recent observation that people’s risk perception, and thus their response, changes over the course of the pandemic. Here, risk response refers to change in human behavior with respect to adherence to NPIs and is one of the sources of pandemic waves. To understand potential approaches for mathematical modeling of human risk response, we reviewed epidemic models that incorporate changes in human response. We narrowed down our search to SIR-type compartmental models, a common and conventional style of modeling in epidemiology. 

\subsection{Findings and Contributions}

During the early stages of our investigation, two exceedingly different approaches to modeling risk response were recognizable: exogenous versus endogenous. In the exogenous approach, models have parameters that incorporate human behavior, such as contact rate, infectivity, or the reproductive number, which may change over the period of an epidemic. The values of the parameters are inputs to the model. On the other hand, the endogenous approach refers to models that formulate human response as a function of a disease state. In this approach, parameter functions change internally due to changes in the disease state, such as the number of cases altering contact rate, representing how society responds to evolving risks. This distinction has major implications for model performance.

Historically, SIR-type models consider human behavior exogenously, where risk response does not depend upon changes internal to the model. Since individual responses are not modeled, modelers have to assume specific scenarios for human behavior. This is shown to significantly affect policy insights as a policy can influence the state of a disease and thus risk perception and people's adherence to NPIs. Furthermore, projecting the state of the disease requires projecting human compliance with NPIs \cite{rahmandad2022enhancing}.  The demand for more modeling of human behavior, driven by the recent pandemic, led to modelers incorporating risk response endogenously. We noted that various approaches are used for endogenous formulation of human response, and these formulations can be further differentiated into two major categories (implicit information diffusion versus explicit information diffusion) based on how information diffusion about the state of the disease is modeled. With implicit information diffusion, change in risk response is dependent upon the state of the disease (or its lagged values) while with explicit information diffusion, more details are incorporated on how people become aware of the risks and spread the information.  

Our review suggests that endogenous models have more capabilities and less limitations than exogenous models with respect to reproducing realistic epidemic behavior, consistent with findings from the systematic review of \citeA{hamilton2024incorporating}. 
Endogenous models are capable of producing multiple epidemic waves, such as those observed during the COVID-19 pandemic and a multitude of other epidemics. A main driving factor behind this capability to produce realistic waves is the negative feedback loop between human behavior and disease dynamics. In exogenous models, any changes in transmission rate are driven externally, limiting their ability to produce multiple epidemic waves without additional complexities such as waning immunity. Furthermore, endogenous construction of risk response allows for the consideration of time delays between event occurrence and public awareness of events (such as infections or deaths), more accurately capturing the disparity between risk perception and actual risk. Simulations show that these time delays further drive oscillatory behavior.

This research makes twofold contributions to the literature. First, by responding to the increasing attention on incorporating human behavior \cite{gentili2020challenges, sooknanan2023fomo, funk2010modelling, perra2021non}, it contributes to the study of epidemic modeling. Past studies provided conceptual frameworks for incorporating human behavior and discussed the importance of coupling human behavior dynamics with epidemic dynamics \cite{ferguson2007capturing, bauch2013social}. This paper adds to the literature by articulating how such frameworks can be formalized and by elaborating the differences between various approaches to modeling changes in human behavior. Furthermore, it shows how the behavioral side can be modeled through implicit or explicit information diffusion formulations. Second, this research contributes to the literature of operations research by examining endogenous modeling approaches in epidemiological models. Many operations research scientists in the past have taken an endogenous modeling approach in modeling various social systems (e.g., \citeA{petropoulos2023operational, lane1998greater}). Particularly in the field of system dynamics, several scholars in the past have emphasized the importance of endogenous modeling, especially in representing human behavior \cite{richardson2011reflections, forrester1997industrial}. In this sense, this study resonates with that particular literature, offering novel evidence from epidemic modeling. Furthermore, by improving epidemic modeling, this study has the potential to influence various optimization studies that depend on epidemic forecasts (e.g, \citeA{kraft2023pandemic, fattahi2023resource, aljuneidi2024forecasting}). For instance, medical supply chain studies aiming to optimize medical resource allocation during epidemics (e.g., \citeA{hammami2023government}) can benefit from more precise projections offered by endogenous behavioral models.

\subsection{Future Studies}

The focus of this literature review is on the human risk-response feedback loop, which posits that as the prevalence of a disease increases, society collectively perceives a higher level of risk. This heightened risk perception leads to activities that decrease the spread of the disease, primarily through NPIs. It is important to note that while risk response is the most commonly studied behavioral mechanism, it is not the only behavioral feedback loop \cite{gordondeveloping2024}. Here, we offer several examples of other possible feedback loops noted in our literature review for future exploration.

A similar behavioral mechanism to risk response is the effect of risk perception on vaccination. Past studies indicated that as risks rise, people are more likely to get vaccinated. It is argued that an individual's experience with past infectious diseases, such as the last seasonal flu, can increase the likelihood of vaccination. This reaction was incorporated in a few studies \cite{yan2023coupling}. 
A major point to consider in formulating willingness to vaccinate is that, in addition to perceived risks of the current epidemic and past epidemic experiences, one's willingness to vaccinate is also affected by available information about the vaccines, their effectiveness, and their possible side effects. Modeling the diffusion of information (and misinformation) about vaccines can also add another layer of complexity. Furthermore, vaccinated individuals may feel safe and increase their social interactions. In a realistic setting where vaccines are not fully effective, this behavior can contribute to new cases.

Another important behavioral mechanism to model is adherence fatigue, the idea that NPI compliance is psychologically and economically costly, and can lead to a decrease in trust in government policies and ultimately disregarding the risks of the pandemic \cite{jorgensen2022pandemic}. This feedback loop can potentially cause large surges of an infectious disease after a prolonged period of compliance. An important aspect of this mechanism is its varying effect on different communities with different economic or psychological tolerance, possibly disproportionately penalizing economically under-served communities. Based on our review, we found a few studies that modeled adherence fatigue \cite{rahmandad2021behavioral,pant2024mathematical}. Other studies found that public trust in government policies can decline over the course of the pandemic leading to less NPI compliance \cite{jorgensen2022pandemic}. Overall, mathematical modeling of trust dynamics during a pandemic is an interesting and potentially promising challenge with major implications for policy implementation and adherence to public health recommendations.

Other behavioral mechanisms to include in modeling encompass the effect of learning, both at the public level and within healthcare systems. Often, early on, there are many unknowns about a disease, leading to more restrictive policies such as lock-downs. Once additional information is acquired about a particular disease, societal reactions to risks move toward less costly and more effective approaches. For example, during the COVID-19 pandemic, after a year, mobility was substantially restored, and people used masking or better air ventilation to reduce the spread of the disease. Additionally, the healthcare system learned how to better treat patients, decreasing the infection fatality rate. Such achievements not only decreased the death rate but also affected the described risk-response feedback loop since the risks declined \cite{rahmandad2021behavioral}.

The described behavioral mechanisms, while often absent in most models, are worth examining. They are not necessarily all active during the entire period of a pandemic. For example, misinformation about vaccines primarily arises after vaccines become available (e.g., see \citeA{mumtaz2022exploring}), or adherence fatigue may occur after the first few waves of the disease (e.g., see \citeA{rahmandad2021behavioral}). Incorporating such feedback loops depends on the purpose of modeling, period of analysis, and the specific problem at hand. 

\subsection{Limitations}

This study has several limitations. First, this work is not intended to present an exhaustive review, as the scope of literature is too vast, but rather to present classifications from the literature which exemplify the importance of endogenous formulations of human behavior of risk response in epidemic disease modeling. Additionally, we focus on a specific mathematical structure of compartmental SIR-type models, although many other model types, e.g., agent-based, network, machine learning, and several more exist \cite{hamilton2024incorporating}. Particularly, an analysis of agent-based models regarding the incorporation of risk response can be informative. In such contexts, modelers will be able to more extensively capture heterogeneities in risk-avoidance behavior, high spreaders' behavior and response to risk, and the potential impacts of individuals on the spread of the disease.

Finally, we limit our examination to models which include human behavior in the form of risk response. Human behavior is incredibly complex and many aspects can be considered, including mobility, social interactions, compliance with preventative measures, willingness to vaccinate, adherence fatigue, social and political interactions and leanings, misinformation, and others. A continued investigation of the various aspects of human behavior within disparate types of mathematical models is critical to continuing to develop understanding of disease spread and allow for better preparedness for the next pandemic.

\subsection{Conclusion}

Overall, this study articulated different approaches to and the usefulness of modeling human behavior in conjunction with human disease models with a focus on societal response to risks. Particularly, an endogenous approach can improve model adaptation to human behaviors changing with respect to the state of the disease (and vice versa). As noted from the evidence in the literature \cite{rahmandad2022enhancing,weitz2020awareness}, the endogenous formulation of risk response can potentially lead to more accurate models and, hence, better long-term projections and improved policy suggestions. A simple approach to model risk response was discussed and various simulation results were depicted.

Our findings have the following implications for modelers. The consideration of human behavior is critical for accurate simulation of epidemics as there is inherent feedback between disease dynamics and human behavior. As such, linking human risk response with disease dynamics, internal to the model, allows for a wider range of model capabilities, such as reproduction of multiple waves and improved forecasting. Furthermore, models do not necessarily need to increase in complexity in order to more closely reflect past outbreaks and offer informative forecasts. Simple models, with the relevant assumptions and formulation, can provide valuable insights of epidemic dynamics.

\section*{Acknowledgements} 
Funding: This work was supported by the US National Science Foundation, Division of Mathematical Sciences and Division of Social and Economic Sciences [2229819].

\section*{Declaration of interest}
Declaration of interest: none.

\section*{Appendix}
\vspace{-1in}
\begin{sidewaystable}[ht!]
\small
\centering
\caption{Parameter values used in simulations for Figures \ref{fig:exog} \& \ref{fig:endo}.} 
\begin{tabular}{@{}clcc@{}}
\hline
\textbf{Figure} & \textbf{Parameter/Description} & \textbf{Value} & \textbf{Source} \\
\hline
6, 7 & $\beta_0$, infectivity rate &  0.7 & \cite{rahmandad2022enhancing} \\
6, 7 & $N$, total population size & $10^6$ & * \\
6, 7 & $\tau_E$, time from infection to infectiousness & 5 & \cite{lauer2020incubation}\\
6, 7 & $\tau_I$, time from infectiousness onset to recovery & 10 & \cite{Centers} \\
6b, 6c & $k$, proportionality constant &  3/500 & *\\
7a, 7b, 7c & $\tau_F$, delay between infection occurrence and information dissemination & 20 & \cite{rahmandad2022enhancing}\\
7a & $\gamma$, diminishing impact of infection risk & 2 & \cite{rahmandad2022enhancing}\\
7a & $\alpha$, community sensitivity to prevalence & 1/30 & \cite{rahmandad2022enhancing}\\
7b & $k$, proportionality constant &  0.05 & *\\
7c & $c$, cost of preventative measures & 0.4 & \cite{Martcheva}\\
7c & $m$, infectivity multiplier & 0.5 & \cite{Martcheva}\\
\hline
\end{tabular}
\label{table:params}
\caption*{* values chosen for visible epidemic curves in a reasonable population size}
\end{sidewaystable}

\clearpage
\newpage
\bibliographystyle{apacite}
\bibliography{References}

\begin{thebibliography}{}

\bibitem [\protect \citeauthoryear {%
Abbas%
, MA%
, Park%
, Parveen%
\BCBL {}\ \BBA {} Kim%
}{%
Abbas%
\ \protect \BOthers {.}}{%
{\protect \APACyear {2022}}%
}]{%
abbas2022evolution}
\APACinsertmetastar {%
abbas2022evolution}%
\begin{APACrefauthors}%
Abbas, W.%
, MA, M.%
, Park, A.%
, Parveen, S.%
\BCBL {}\ \BBA {} Kim, S.%
\end{APACrefauthors}%
\unskip\
\newblock
\APACrefYearMonthDay{2022}{}{}.
\newblock
{\BBOQ}\APACrefatitle {Evolution and consequences of individual responses
  during the {COVID-19 outbreak}} {Evolution and consequences of individual
  responses during the {COVID-19 outbreak}}.{\BBCQ}
\newblock
\APACjournalVolNumPages{PLOS One}{17}{9}{e0273964}.
\PrintBackRefs{\CurrentBib}

\bibitem [\protect \citeauthoryear {%
Alharbi%
\ \BBA {} Kribs%
}{%
Alharbi%
\ \BBA {} Kribs%
}{%
{\protect \APACyear {2022}}%
}]{%
alharbi2022nature}
\APACinsertmetastar {%
alharbi2022nature}%
\begin{APACrefauthors}%
Alharbi, M\BPBI H.%
\BCBT {}\ \BBA {} Kribs, C\BPBI M.%
\end{APACrefauthors}%
\unskip\
\newblock
\APACrefYearMonthDay{2022}{}{}.
\newblock
{\BBOQ}\APACrefatitle {How the nature of behavior change affects the impact of
  asymptomatic coronavirus transmission} {How the nature of behavior change
  affects the impact of asymptomatic coronavirus transmission}.{\BBCQ}
\newblock
\APACjournalVolNumPages{Ricerche di Matematica}{}{}{1--21}.
\PrintBackRefs{\CurrentBib}

\bibitem [\protect \citeauthoryear {%
Aljuneidi%
, Punia%
, Jebali%
\BCBL {}\ \BBA {} Nikolopoulos%
}{%
Aljuneidi%
\ \protect \BOthers {.}}{%
{\protect \APACyear {2024}}%
}]{%
aljuneidi2024forecasting}
\APACinsertmetastar {%
aljuneidi2024forecasting}%
\begin{APACrefauthors}%
Aljuneidi, T.%
, Punia, S.%
, Jebali, A.%
\BCBL {}\ \BBA {} Nikolopoulos, K.%
\end{APACrefauthors}%
\unskip\
\newblock
\APACrefYearMonthDay{2024}{}{}.
\newblock
{\BBOQ}\APACrefatitle {Forecasting and planning for a critical infrastructure
  sector during a pandemic: Empirical evidence from a food supply chain}
  {Forecasting and planning for a critical infrastructure sector during a
  pandemic: Empirical evidence from a food supply chain}.{\BBCQ}
\newblock
\APACjournalVolNumPages{European Journal of Operational
  Research}{317}{3}{936--952}.
\PrintBackRefs{\CurrentBib}

\bibitem [\protect \citeauthoryear {%
Anderson%
\ \BBA {} May%
}{%
Anderson%
\ \BBA {} May%
}{%
{\protect \APACyear {1979}}%
}]{%
anderson1979population}
\APACinsertmetastar {%
anderson1979population}%
\begin{APACrefauthors}%
Anderson, R\BPBI M.%
\BCBT {}\ \BBA {} May, R\BPBI M.%
\end{APACrefauthors}%
\unskip\
\newblock
\APACrefYearMonthDay{1979}{}{}.
\newblock
{\BBOQ}\APACrefatitle {Population biology of infectious diseases: Part I}
  {Population biology of infectious diseases: Part i}.{\BBCQ}
\newblock
\APACjournalVolNumPages{Nature}{280}{5721}{361--367}.
\PrintBackRefs{\CurrentBib}

\bibitem [\protect \citeauthoryear {%
Anderson%
\ \BBA {} May%
}{%
Anderson%
\ \BBA {} May%
}{%
{\protect \APACyear {1991}}%
}]{%
anderson1991infectious}
\APACinsertmetastar {%
anderson1991infectious}%
\begin{APACrefauthors}%
Anderson, R\BPBI M.%
\BCBT {}\ \BBA {} May, R\BPBI M.%
\end{APACrefauthors}%
\unskip\
\newblock
\APACrefYear{1991}.
\newblock
\APACrefbtitle {Infectious diseases of humans: dynamics and control}
  {Infectious diseases of humans: dynamics and control}.
\newblock
\APACaddressPublisher{}{Oxford University Press}.
\PrintBackRefs{\CurrentBib}

\bibitem [\protect \citeauthoryear {%
Bass%
}{%
Bass%
}{%
{\protect \APACyear {1969}}%
}]{%
bass1969new}
\APACinsertmetastar {%
bass1969new}%
\begin{APACrefauthors}%
Bass, F\BPBI M.%
\end{APACrefauthors}%
\unskip\
\newblock
\APACrefYearMonthDay{1969}{}{}.
\newblock
{\BBOQ}\APACrefatitle {A new product growth for model consumer durables} {A new
  product growth for model consumer durables}.{\BBCQ}
\newblock
\APACjournalVolNumPages{Management Science}{15}{5}{215--227}.
\PrintBackRefs{\CurrentBib}

\bibitem [\protect \citeauthoryear {%
Bauch%
\ \BBA {} Galvani%
}{%
Bauch%
\ \BBA {} Galvani%
}{%
{\protect \APACyear {2013}}%
}]{%
bauch2013social}
\APACinsertmetastar {%
bauch2013social}%
\begin{APACrefauthors}%
Bauch, C\BPBI T.%
\BCBT {}\ \BBA {} Galvani, A\BPBI P.%
\end{APACrefauthors}%
\unskip\
\newblock
\APACrefYearMonthDay{2013}{}{}.
\newblock
{\BBOQ}\APACrefatitle {Social factors in epidemiology} {Social factors in
  epidemiology}.{\BBCQ}
\newblock
\APACjournalVolNumPages{Science}{342}{6154}{47--49}.
\PrintBackRefs{\CurrentBib}

\bibitem [\protect \citeauthoryear {%
Bliman%
, Carrozzo-Magli%
, D'onofrio%
\BCBL {}\ \BBA {} Manfredi%
}{%
Bliman%
\ \protect \BOthers {.}}{%
{\protect \APACyear {2022}}%
}]{%
bliman2022tiered}
\APACinsertmetastar {%
bliman2022tiered}%
\begin{APACrefauthors}%
Bliman, P\BHBI A.%
, Carrozzo-Magli, A.%
, D'onofrio, A.%
\BCBL {}\ \BBA {} Manfredi, P.%
\end{APACrefauthors}%
\unskip\
\newblock
\APACrefYearMonthDay{2022}{}{}.
\newblock
{\BBOQ}\APACrefatitle {Tiered social distancing policies and epidemic control}
  {Tiered social distancing policies and epidemic control}.{\BBCQ}
\newblock
\APACjournalVolNumPages{Proceedings of the Royal Society
  A}{478}{2268}{20220175}.
\PrintBackRefs{\CurrentBib}

\bibitem [\protect \citeauthoryear {%
Browne%
, Gulbudak%
\BCBL {}\ \BBA {} Macdonald%
}{%
Browne%
\ \protect \BOthers {.}}{%
{\protect \APACyear {2022}}%
}]{%
browne2022differential}
\APACinsertmetastar {%
browne2022differential}%
\begin{APACrefauthors}%
Browne, C\BPBI J.%
, Gulbudak, H.%
\BCBL {}\ \BBA {} Macdonald, J\BPBI C.%
\end{APACrefauthors}%
\unskip\
\newblock
\APACrefYearMonthDay{2022}{}{}.
\newblock
{\BBOQ}\APACrefatitle {Differential impacts of contact tracing and lockdowns on
  outbreak size in {COVID-19 model applied to China}} {Differential impacts of
  contact tracing and lockdowns on outbreak size in {COVID-19 model applied to
  China}}.{\BBCQ}
\newblock
\APACjournalVolNumPages{Journal of Theoretical Biology}{532}{}{110919}.
\PrintBackRefs{\CurrentBib}

\bibitem [\protect \citeauthoryear {%
Bubar%
\ \protect \BOthers {.}}{%
Bubar%
\ \protect \BOthers {.}}{%
{\protect \APACyear {2021}}%
}]{%
bubar2021model}
\APACinsertmetastar {%
bubar2021model}%
\begin{APACrefauthors}%
Bubar, K\BPBI M.%
, Reinholt, K.%
, Kissler, S\BPBI M.%
, Lipsitch, M.%
, Cobey, S.%
, Grad, Y\BPBI H.%
\BCBL {}\ \BBA {} Larremore, D\BPBI B.%
\end{APACrefauthors}%
\unskip\
\newblock
\APACrefYearMonthDay{2021}{}{}.
\newblock
{\BBOQ}\APACrefatitle {Model-informed COVID-19 vaccine prioritization
  strategies by age and serostatus} {Model-informed covid-19 vaccine
  prioritization strategies by age and serostatus}.{\BBCQ}
\newblock
\APACjournalVolNumPages{Science}{371}{6352}{916--921}.
\PrintBackRefs{\CurrentBib}

\bibitem [\protect \citeauthoryear {%
Buonomo%
\ \BBA {} Della~Marca%
}{%
Buonomo%
\ \BBA {} Della~Marca%
}{%
{\protect \APACyear {2020}}%
}]{%
buonomo2020effects}
\APACinsertmetastar {%
buonomo2020effects}%
\begin{APACrefauthors}%
Buonomo, B.%
\BCBT {}\ \BBA {} Della~Marca, R.%
\end{APACrefauthors}%
\unskip\
\newblock
\APACrefYearMonthDay{2020}{}{}.
\newblock
{\BBOQ}\APACrefatitle {Effects of information-induced behavioural changes
  during the {COVID-19 lockdowns: the case of Italy}} {Effects of
  information-induced behavioural changes during the {COVID-19 lockdowns: the
  case of Italy}}.{\BBCQ}
\newblock
\APACjournalVolNumPages{Royal Society Open Science}{7}{10}{201635}.
\PrintBackRefs{\CurrentBib}

\bibitem [\protect \citeauthoryear {%
\APACcitebtitle {Centers for Disease Control and Prevention 2023}}{%
\APACcitebtitle {Centers for Disease Control and Prevention 2023}}{%
{\protect \APACyear {2023}}%
}]{%
Centers}
\APACinsertmetastar {%
Centers}%
\APACrefbtitle {Centers for Disease Control and Prevention 2023.} {Centers for
  disease control and prevention 2023.}
\newblock
\APACrefYearMonthDay{2023}{Aug}{}.
\newblock
\APACaddressPublisher{}{Centers for Disease Control and Prevention}.
\newblock
\begin{APACrefURL}
  [{2024-02-05}]\url{https://www.cdc.gov/coronavirus/2019-ncov/hcp/duration-isolation.html}
  \end{APACrefURL}
\PrintBackRefs{\CurrentBib}

\bibitem [\protect \citeauthoryear {%
Childs%
\ \protect \BOthers {.}}{%
Childs%
\ \protect \BOthers {.}}{%
{\protect \APACyear {2022}}%
}]{%
childs2022modeling}
\APACinsertmetastar {%
childs2022modeling}%
\begin{APACrefauthors}%
Childs, L.%
, Dick, D\BPBI W.%
, Feng, Z.%
, Heffernan, J\BPBI M.%
, Li, J.%
\BCBL {}\ \BBA {} R{\"o}st, G.%
\end{APACrefauthors}%
\unskip\
\newblock
\APACrefYearMonthDay{2022}{}{}.
\newblock
{\BBOQ}\APACrefatitle {Modeling waning and boosting of COVID-19 in Canada with
  vaccination} {Modeling waning and boosting of covid-19 in canada with
  vaccination}.{\BBCQ}
\newblock
\APACjournalVolNumPages{Epidemics}{39}{}{100583}.
\PrintBackRefs{\CurrentBib}

\bibitem [\protect \citeauthoryear {%
Costello%
, Watts%
\BCBL {}\ \BBA {} Howe%
}{%
Costello%
\ \protect \BOthers {.}}{%
{\protect \APACyear {2023}}%
}]{%
costello2023model}
\APACinsertmetastar {%
costello2023model}%
\begin{APACrefauthors}%
Costello, F.%
, Watts, P.%
\BCBL {}\ \BBA {} Howe, R.%
\end{APACrefauthors}%
\unskip\
\newblock
\APACrefYearMonthDay{2023}{}{}.
\newblock
{\BBOQ}\APACrefatitle {A model of behavioural response to risk accurately
  predicts the statistical distribution of {COVID-19} infection and
  reproduction numbers} {A model of behavioural response to risk accurately
  predicts the statistical distribution of {COVID-19} infection and
  reproduction numbers}.{\BBCQ}
\newblock
\APACjournalVolNumPages{Scientific Reports}{13}{1}{2435}.
\PrintBackRefs{\CurrentBib}

\bibitem [\protect \citeauthoryear {%
Currie%
\ \protect \BOthers {.}}{%
Currie%
\ \protect \BOthers {.}}{%
{\protect \APACyear {2020}}%
}]{%
currie2020simulation}
\APACinsertmetastar {%
currie2020simulation}%
\begin{APACrefauthors}%
Currie, C\BPBI S.%
, Fowler, J\BPBI W.%
, Kotiadis, K.%
, Monks, T.%
, Onggo, B\BPBI S.%
, Robertson, D\BPBI A.%
\BCBL {}\ \BBA {} Tako, A\BPBI A.%
\end{APACrefauthors}%
\unskip\
\newblock
\APACrefYearMonthDay{2020}{}{}.
\newblock
{\BBOQ}\APACrefatitle {How simulation modelling can help reduce the impact of
  COVID-19} {How simulation modelling can help reduce the impact of
  covid-19}.{\BBCQ}
\newblock
\APACjournalVolNumPages{Journal of Simulation}{14}{2}{83--97}.
\PrintBackRefs{\CurrentBib}

\bibitem [\protect \citeauthoryear {%
De~Meijere%
, Colizza%
, Valdano%
\BCBL {}\ \BBA {} Castellano%
}{%
De~Meijere%
\ \protect \BOthers {.}}{%
{\protect \APACyear {2021}}%
}]{%
de2021effect}
\APACinsertmetastar {%
de2021effect}%
\begin{APACrefauthors}%
De~Meijere, G.%
, Colizza, V.%
, Valdano, E.%
\BCBL {}\ \BBA {} Castellano, C.%
\end{APACrefauthors}%
\unskip\
\newblock
\APACrefYearMonthDay{2021}{}{}.
\newblock
{\BBOQ}\APACrefatitle {Effect of delayed awareness and fatigue on the efficacy
  of self-isolation in epidemic control} {Effect of delayed awareness and
  fatigue on the efficacy of self-isolation in epidemic control}.{\BBCQ}
\newblock
\APACjournalVolNumPages{Physical Review E}{104}{4}{044316}.
\PrintBackRefs{\CurrentBib}

\bibitem [\protect \citeauthoryear {%
Dick%
\ \protect \BOthers {.}}{%
Dick%
\ \protect \BOthers {.}}{%
{\protect \APACyear {2021}}%
}]{%
dick2021covid}
\APACinsertmetastar {%
dick2021covid}%
\begin{APACrefauthors}%
Dick, D\BPBI W.%
, Childs, L.%
, Feng, Z.%
, Li, J.%
, R{\"o}st, G.%
, Buckeridge, D\BPBI L.%
\BDBL {}Heffernan, J\BPBI M.%
\end{APACrefauthors}%
\unskip\
\newblock
\APACrefYearMonthDay{2021}{}{}.
\newblock
{\BBOQ}\APACrefatitle {COVID-19 seroprevalence in Canada modelling waning and
  boosting COVID-19 immunity in Canada a Canadian immunization research network
  study} {Covid-19 seroprevalence in canada modelling waning and boosting
  covid-19 immunity in canada a canadian immunization research network
  study}.{\BBCQ}
\newblock
\APACjournalVolNumPages{Vaccines}{10}{1}{17}.
\PrintBackRefs{\CurrentBib}

\bibitem [\protect \citeauthoryear {%
d'Onofrio%
\ \BBA {} Manfredi%
}{%
d'Onofrio%
\ \BBA {} Manfredi%
}{%
{\protect \APACyear {2022}}%
}]{%
d2022behavioral}
\APACinsertmetastar {%
d2022behavioral}%
\begin{APACrefauthors}%
d'Onofrio, A.%
\BCBT {}\ \BBA {} Manfredi, P.%
\end{APACrefauthors}%
\unskip\
\newblock
\APACrefYearMonthDay{2022}{}{}.
\newblock
{\BBOQ}\APACrefatitle {Behavioral {SIR} models with incidence-based
  social-distancing} {Behavioral {SIR} models with incidence-based
  social-distancing}.{\BBCQ}
\newblock
\APACjournalVolNumPages{Chaos, Solitons \& Fractals}{159}{}{112072}.
\PrintBackRefs{\CurrentBib}

\bibitem [\protect \citeauthoryear {%
Duggan%
\ \protect \BOthers {.}}{%
Duggan%
\ \protect \BOthers {.}}{%
{\protect \APACyear {2024}}%
}]{%
duggan2024age}
\APACinsertmetastar {%
duggan2024age}%
\begin{APACrefauthors}%
Duggan, J.%
, Andrade, J.%
, Murphy, T\BPBI B.%
, Gleeson, J\BPBI P.%
, Walsh, C.%
\BCBL {}\ \BBA {} Nolan, P.%
\end{APACrefauthors}%
\unskip\
\newblock
\APACrefYearMonthDay{2024}{}{}.
\newblock
{\BBOQ}\APACrefatitle {An age-cohort simulation model for generating COVID-19
  scenarios: A study from Ireland's pandemic response} {An age-cohort
  simulation model for generating covid-19 scenarios: A study from ireland's
  pandemic response}.{\BBCQ}
\newblock
\APACjournalVolNumPages{European Journal of Operational
  Research}{313}{1}{343--358}.
\PrintBackRefs{\CurrentBib}

\bibitem [\protect \citeauthoryear {%
Epstein%
, Hatna%
\BCBL {}\ \BBA {} Crodelle%
}{%
Epstein%
\ \protect \BOthers {.}}{%
{\protect \APACyear {2021}}%
}]{%
epstein2021triple}
\APACinsertmetastar {%
epstein2021triple}%
\begin{APACrefauthors}%
Epstein, J\BPBI M.%
, Hatna, E.%
\BCBL {}\ \BBA {} Crodelle, J.%
\end{APACrefauthors}%
\unskip\
\newblock
\APACrefYearMonthDay{2021}{}{}.
\newblock
{\BBOQ}\APACrefatitle {Triple contagion: a two-fears epidemic model} {Triple
  contagion: a two-fears epidemic model}.{\BBCQ}
\newblock
\APACjournalVolNumPages{Journal of the Royal Society
  Interface}{18}{181}{20210186}.
\PrintBackRefs{\CurrentBib}

\bibitem [\protect \citeauthoryear {%
Eryarsoy%
, Shahmanzari%
\BCBL {}\ \BBA {} Tanrisever%
}{%
Eryarsoy%
\ \protect \BOthers {.}}{%
{\protect \APACyear {2023}}%
}]{%
eryarsoy2023models}
\APACinsertmetastar {%
eryarsoy2023models}%
\begin{APACrefauthors}%
Eryarsoy, E.%
, Shahmanzari, M.%
\BCBL {}\ \BBA {} Tanrisever, F.%
\end{APACrefauthors}%
\unskip\
\newblock
\APACrefYearMonthDay{2023}{}{}.
\newblock
{\BBOQ}\APACrefatitle {Models for government intervention during a pandemic}
  {Models for government intervention during a pandemic}.{\BBCQ}
\newblock
\APACjournalVolNumPages{European Journal of Operational
  Research}{304}{1}{69--83}.
\PrintBackRefs{\CurrentBib}

\bibitem [\protect \citeauthoryear {%
Fattahi%
, Keyvanshokooh%
, Kannan%
\BCBL {}\ \BBA {} Govindan%
}{%
Fattahi%
\ \protect \BOthers {.}}{%
{\protect \APACyear {2023}}%
}]{%
fattahi2023resource}
\APACinsertmetastar {%
fattahi2023resource}%
\begin{APACrefauthors}%
Fattahi, M.%
, Keyvanshokooh, E.%
, Kannan, D.%
\BCBL {}\ \BBA {} Govindan, K.%
\end{APACrefauthors}%
\unskip\
\newblock
\APACrefYearMonthDay{2023}{}{}.
\newblock
{\BBOQ}\APACrefatitle {Resource planning strategies for healthcare systems
  during a pandemic} {Resource planning strategies for healthcare systems
  during a pandemic}.{\BBCQ}
\newblock
\APACjournalVolNumPages{European Journal of Operational
  Research}{304}{1}{192--206}.
\PrintBackRefs{\CurrentBib}

\bibitem [\protect \citeauthoryear {%
Ferguson%
}{%
Ferguson%
}{%
{\protect \APACyear {2007}}%
}]{%
ferguson2007capturing}
\APACinsertmetastar {%
ferguson2007capturing}%
\begin{APACrefauthors}%
Ferguson, N.%
\end{APACrefauthors}%
\unskip\
\newblock
\APACrefYearMonthDay{2007}{}{}.
\newblock
{\BBOQ}\APACrefatitle {Capturing human behaviour} {Capturing human
  behaviour}.{\BBCQ}
\newblock
\APACjournalVolNumPages{Nature}{446}{7137}{733--733}.
\PrintBackRefs{\CurrentBib}

\bibitem [\protect \citeauthoryear {%
Forrester%
}{%
Forrester%
}{%
{\protect \APACyear {1997}}%
}]{%
forrester1997industrial}
\APACinsertmetastar {%
forrester1997industrial}%
\begin{APACrefauthors}%
Forrester, J\BPBI W.%
\end{APACrefauthors}%
\unskip\
\newblock
\APACrefYearMonthDay{1997}{}{}.
\newblock
{\BBOQ}\APACrefatitle {Industrial dynamics} {Industrial dynamics}.{\BBCQ}
\newblock
\APACjournalVolNumPages{Journal of the Operational Research
  Society}{48}{10}{1037--1041}.
\PrintBackRefs{\CurrentBib}

\bibitem [\protect \citeauthoryear {%
Funk%
, Salath{\'e}%
\BCBL {}\ \BBA {} Jansen%
}{%
Funk%
\ \protect \BOthers {.}}{%
{\protect \APACyear {2010}}%
}]{%
funk2010modelling}
\APACinsertmetastar {%
funk2010modelling}%
\begin{APACrefauthors}%
Funk, S.%
, Salath{\'e}, M.%
\BCBL {}\ \BBA {} Jansen, V\BPBI A.%
\end{APACrefauthors}%
\unskip\
\newblock
\APACrefYearMonthDay{2010}{}{}.
\newblock
{\BBOQ}\APACrefatitle {Modelling the influence of human behaviour on the spread
  of infectious diseases: a review} {Modelling the influence of human behaviour
  on the spread of infectious diseases: a review}.{\BBCQ}
\newblock
\APACjournalVolNumPages{Journal of the Royal Society
  Interface}{7}{50}{1247--1256}.
\PrintBackRefs{\CurrentBib}

\bibitem [\protect \citeauthoryear {%
Gentili%
\ \BBA {} Cristea%
}{%
Gentili%
\ \BBA {} Cristea%
}{%
{\protect \APACyear {2020}}%
}]{%
gentili2020challenges}
\APACinsertmetastar {%
gentili2020challenges}%
\begin{APACrefauthors}%
Gentili, C.%
\BCBT {}\ \BBA {} Cristea, I\BPBI A.%
\end{APACrefauthors}%
\unskip\
\newblock
\APACrefYearMonthDay{2020}{}{}.
\newblock
{\BBOQ}\APACrefatitle {Challenges and opportunities for human behavior research
  in the coronavirus disease (COVID-19) pandemic} {Challenges and opportunities
  for human behavior research in the coronavirus disease (covid-19)
  pandemic}.{\BBCQ}
\newblock
\APACjournalVolNumPages{Frontiers in psychology}{11}{}{1786}.
\PrintBackRefs{\CurrentBib}

\bibitem [\protect \citeauthoryear {%
Ghaffarzadegan%
}{%
Ghaffarzadegan%
}{%
{\protect \APACyear {2021}}%
}]{%
ghaffarzadegan2021simulation}
\APACinsertmetastar {%
ghaffarzadegan2021simulation}%
\begin{APACrefauthors}%
Ghaffarzadegan, N.%
\end{APACrefauthors}%
\unskip\
\newblock
\APACrefYearMonthDay{2021}{}{}.
\newblock
{\BBOQ}\APACrefatitle {Simulation-based what-if analysis for controlling the
  spread of {COVID-19} in universities} {Simulation-based what-if analysis for
  controlling the spread of {COVID-19} in universities}.{\BBCQ}
\newblock
\APACjournalVolNumPages{PLOS One}{16}{2}{e0246323}.
\PrintBackRefs{\CurrentBib}

\bibitem [\protect \citeauthoryear {%
Gordon%
\ \protect \BOthers {.}}{%
Gordon%
\ \protect \BOthers {.}}{%
{\protect \APACyear {2024}}%
}]{%
gordondeveloping2024}
\APACinsertmetastar {%
gordondeveloping2024}%
\begin{APACrefauthors}%
Gordon, D.%
, Mashayekhi, A\BPBI N.%
, Tomoaia-Cotisel, A.%
, Kim, H.%
, Bahaddin, B.%
, Luna-Reyes, L\BPBI F.%
\BCBL {}\ \BBA {} Andersen, D\BPBI F.%
\end{APACrefauthors}%
\unskip\
\newblock
\APACrefYearMonthDay{2024}{}{}.
\newblock
{\BBOQ}\APACrefatitle {Developing model-based storytelling to share systemic
  insights to the public during the {COVID-19} pandemic} {Developing
  model-based storytelling to share systemic insights to the public during the
  {COVID-19} pandemic}.{\BBCQ}
\newblock
\APACjournalVolNumPages{System Dynamics Review}{}{}{}.
\PrintBackRefs{\CurrentBib}

\bibitem [\protect \citeauthoryear {%
Hamilton%
\ \protect \BOthers {.}}{%
Hamilton%
\ \protect \BOthers {.}}{%
{\protect \APACyear {2024}}%
}]{%
hamilton2024incorporating}
\APACinsertmetastar {%
hamilton2024incorporating}%
\begin{APACrefauthors}%
Hamilton, A.%
, Haghpanah, F.%
, Tulchinsky, A.%
, Kipshidze, N.%
, Poleon, S.%
, Lin, G.%
\BDBL {}Klein, E.%
\end{APACrefauthors}%
\unskip\
\newblock
\APACrefYearMonthDay{2024}{}{}.
\newblock
{\BBOQ}\APACrefatitle {Incorporating endogenous human behavior in models of
  COVID-19 transmission: A systematic scoping review} {Incorporating endogenous
  human behavior in models of covid-19 transmission: A systematic scoping
  review}.{\BBCQ}
\newblock
\APACjournalVolNumPages{Dialogues in Health}{}{}{100179}.
\PrintBackRefs{\CurrentBib}

\bibitem [\protect \citeauthoryear {%
Hammami%
, Salman%
, Khouja%
, Nouira%
\BCBL {}\ \BBA {} Alaswad%
}{%
Hammami%
\ \protect \BOthers {.}}{%
{\protect \APACyear {2023}}%
}]{%
hammami2023government}
\APACinsertmetastar {%
hammami2023government}%
\begin{APACrefauthors}%
Hammami, R.%
, Salman, S.%
, Khouja, M.%
, Nouira, I.%
\BCBL {}\ \BBA {} Alaswad, S.%
\end{APACrefauthors}%
\unskip\
\newblock
\APACrefYearMonthDay{2023}{}{}.
\newblock
{\BBOQ}\APACrefatitle {Government strategies to secure the supply of medical
  products in pandemic times} {Government strategies to secure the supply of
  medical products in pandemic times}.{\BBCQ}
\newblock
\APACjournalVolNumPages{European Journal of Operational
  Research}{306}{3}{1364--1387}.
\PrintBackRefs{\CurrentBib}

\bibitem [\protect \citeauthoryear {%
{IHME COVID-19 Forecasting Team}%
}{%
{IHME COVID-19 Forecasting Team}%
}{%
{\protect \APACyear {2021}}%
}]{%
ihme2021modeling}
\APACinsertmetastar {%
ihme2021modeling}%
\begin{APACrefauthors}%
{IHME COVID-19 Forecasting Team}.%
\end{APACrefauthors}%
\unskip\
\newblock
\APACrefYearMonthDay{2021}{}{}.
\newblock
{\BBOQ}\APACrefatitle {Modeling {COVID-19} scenarios for the United States}
  {Modeling {COVID-19} scenarios for the united states}.{\BBCQ}
\newblock
\APACjournalVolNumPages{Nature Medicine}{27}{1}{94--105}.
\PrintBackRefs{\CurrentBib}

\bibitem [\protect \citeauthoryear {%
Jasim%
\ \BBA {} AL-Husseiny%
}{%
Jasim%
\ \BBA {} AL-Husseiny%
}{%
{\protect \APACyear {2021}}%
}]{%
jasim2021study}
\APACinsertmetastar {%
jasim2021study}%
\begin{APACrefauthors}%
Jasim, W\BPBI I.%
\BCBT {}\ \BBA {} AL-Husseiny, H\BPBI F.%
\end{APACrefauthors}%
\unskip\
\newblock
\APACrefYearMonthDay{2021}{}{}.
\newblock
{\BBOQ}\APACrefatitle {A study of shigellosis bacteria disease model with
  awareness effects} {A study of shigellosis bacteria disease model with
  awareness effects}.{\BBCQ}
\newblock
\APACjournalVolNumPages{Ibn AL-Haitham Journal For Pure and Applied
  Sciences}{34}{2}{129--143}.
\PrintBackRefs{\CurrentBib}

\bibitem [\protect \citeauthoryear {%
J{\o}rgensen%
, Bor%
, Rasmussen%
, Lindholt%
\BCBL {}\ \BBA {} Petersen%
}{%
J{\o}rgensen%
\ \protect \BOthers {.}}{%
{\protect \APACyear {2022}}%
}]{%
jorgensen2022pandemic}
\APACinsertmetastar {%
jorgensen2022pandemic}%
\begin{APACrefauthors}%
J{\o}rgensen, F.%
, Bor, A.%
, Rasmussen, M\BPBI S.%
, Lindholt, M\BPBI F.%
\BCBL {}\ \BBA {} Petersen, M\BPBI B.%
\end{APACrefauthors}%
\unskip\
\newblock
\APACrefYearMonthDay{2022}{}{}.
\newblock
{\BBOQ}\APACrefatitle {{Pandemic fatigue fueled political discontent during the
  COVID-19 pandemic}} {{Pandemic fatigue fueled political discontent during the
  COVID-19 pandemic}}.{\BBCQ}
\newblock
\APACjournalVolNumPages{Proceedings of the National Academy of
  Sciences}{119}{48}{e2201266119}.
\PrintBackRefs{\CurrentBib}

\bibitem [\protect \citeauthoryear {%
Juher%
, Rojas%
\BCBL {}\ \BBA {} Salda{\~n}a%
}{%
Juher%
\ \protect \BOthers {.}}{%
{\protect \APACyear {2023}}%
}]{%
juher2023saddle}
\APACinsertmetastar {%
juher2023saddle}%
\begin{APACrefauthors}%
Juher, D.%
, Rojas, D.%
\BCBL {}\ \BBA {} Salda{\~n}a, J.%
\end{APACrefauthors}%
\unskip\
\newblock
\APACrefYearMonthDay{2023}{}{}.
\newblock
{\BBOQ}\APACrefatitle {Saddle--node bifurcation of limit cycles in an epidemic
  model with two levels of awareness} {Saddle--node bifurcation of limit cycles
  in an epidemic model with two levels of awareness}.{\BBCQ}
\newblock
\APACjournalVolNumPages{Physica D: Nonlinear Phenomena}{448}{}{133714}.
\PrintBackRefs{\CurrentBib}

\bibitem [\protect \citeauthoryear {%
Kassa%
, Njagarah%
\BCBL {}\ \BBA {} Terefe%
}{%
Kassa%
\ \protect \BOthers {.}}{%
{\protect \APACyear {2020}}%
}]{%
kassa2020analysis}
\APACinsertmetastar {%
kassa2020analysis}%
\begin{APACrefauthors}%
Kassa, S\BPBI M.%
, Njagarah, J\BPBI B.%
\BCBL {}\ \BBA {} Terefe, Y\BPBI A.%
\end{APACrefauthors}%
\unskip\
\newblock
\APACrefYearMonthDay{2020}{}{}.
\newblock
{\BBOQ}\APACrefatitle {Analysis of the mitigation strategies for {COVID-19}:
  from mathematical modelling perspective} {Analysis of the mitigation
  strategies for {COVID-19}: from mathematical modelling perspective}.{\BBCQ}
\newblock
\APACjournalVolNumPages{Chaos, Solitons \& Fractals}{138}{}{109968}.
\PrintBackRefs{\CurrentBib}

\bibitem [\protect \citeauthoryear {%
Keeling%
\ \BBA {} Rohani%
}{%
Keeling%
\ \BBA {} Rohani%
}{%
{\protect \APACyear {2011}}%
}]{%
keeling2011modeling}
\APACinsertmetastar {%
keeling2011modeling}%
\begin{APACrefauthors}%
Keeling, M\BPBI J.%
\BCBT {}\ \BBA {} Rohani, P.%
\end{APACrefauthors}%
\unskip\
\newblock
\APACrefYear{2011}.
\newblock
\APACrefbtitle {Modeling infectious diseases in humans and animals} {Modeling
  infectious diseases in humans and animals}.
\newblock
\APACaddressPublisher{}{Princeton University Press}.
\PrintBackRefs{\CurrentBib}

\bibitem [\protect \citeauthoryear {%
Kermack%
\ \BBA {} McKendrick%
}{%
Kermack%
\ \BBA {} McKendrick%
}{%
{\protect \APACyear {1927}}%
}]{%
kermack1927contribution}
\APACinsertmetastar {%
kermack1927contribution}%
\begin{APACrefauthors}%
Kermack, W\BPBI O.%
\BCBT {}\ \BBA {} McKendrick, A\BPBI G.%
\end{APACrefauthors}%
\unskip\
\newblock
\APACrefYearMonthDay{1927}{}{}.
\newblock
{\BBOQ}\APACrefatitle {A contribution to the mathematical theory of epidemics}
  {A contribution to the mathematical theory of epidemics}.{\BBCQ}
\newblock
\APACjournalVolNumPages{Proceedings of the Royal Society of London. Series A,
  Containing papers of a mathematical and physical
  character}{115}{772}{700--721}.
\PrintBackRefs{\CurrentBib}

\bibitem [\protect \citeauthoryear {%
Kraft%
\ \BBA {} Weiss%
}{%
Kraft%
\ \BBA {} Weiss%
}{%
{\protect \APACyear {2023}}%
}]{%
kraft2023pandemic}
\APACinsertmetastar {%
kraft2023pandemic}%
\begin{APACrefauthors}%
Kraft, H.%
\BCBT {}\ \BBA {} Weiss, F.%
\end{APACrefauthors}%
\unskip\
\newblock
\APACrefYearMonthDay{2023}{}{}.
\newblock
{\BBOQ}\APACrefatitle {Pandemic portfolio choice} {Pandemic portfolio
  choice}.{\BBCQ}
\newblock
\APACjournalVolNumPages{European Journal of Operational
  Research}{305}{1}{451--462}.
\PrintBackRefs{\CurrentBib}

\bibitem [\protect \citeauthoryear {%
Kumar%
, Sharma%
\BCBL {}\ \BBA {} Singh%
}{%
Kumar%
\ \protect \BOthers {.}}{%
{\protect \APACyear {2023}}%
}]{%
kumar2023multiscale}
\APACinsertmetastar {%
kumar2023multiscale}%
\begin{APACrefauthors}%
Kumar, S.%
, Sharma, B.%
\BCBL {}\ \BBA {} Singh, V.%
\end{APACrefauthors}%
\unskip\
\newblock
\APACrefYearMonthDay{2023}{}{}.
\newblock
{\BBOQ}\APACrefatitle {A multiscale modeling framework to study the
  interdependence of brain, behavior, and pandemic} {A multiscale modeling
  framework to study the interdependence of brain, behavior, and
  pandemic}.{\BBCQ}
\newblock
\APACjournalVolNumPages{Nonlinear Dynamics}{111}{8}{7729--7749}.
\PrintBackRefs{\CurrentBib}

\bibitem [\protect \citeauthoryear {%
Lane%
\ \BBA {} Oliva%
}{%
Lane%
\ \BBA {} Oliva%
}{%
{\protect \APACyear {1998}}%
}]{%
lane1998greater}
\APACinsertmetastar {%
lane1998greater}%
\begin{APACrefauthors}%
Lane, D\BPBI C.%
\BCBT {}\ \BBA {} Oliva, R.%
\end{APACrefauthors}%
\unskip\
\newblock
\APACrefYearMonthDay{1998}{}{}.
\newblock
{\BBOQ}\APACrefatitle {The greater whole: Towards a synthesis of system
  dynamics and soft systems methodology} {The greater whole: Towards a
  synthesis of system dynamics and soft systems methodology}.{\BBCQ}
\newblock
\APACjournalVolNumPages{European journal of operational
  research}{107}{1}{214--235}.
\PrintBackRefs{\CurrentBib}

\bibitem [\protect \citeauthoryear {%
Lauer%
\ \protect \BOthers {.}}{%
Lauer%
\ \protect \BOthers {.}}{%
{\protect \APACyear {2020}}%
}]{%
lauer2020incubation}
\APACinsertmetastar {%
lauer2020incubation}%
\begin{APACrefauthors}%
Lauer, S\BPBI A.%
, Grantz, K\BPBI H.%
, Bi, Q.%
, Jones, F\BPBI K.%
, Zheng, Q.%
, Meredith, H\BPBI R.%
\BDBL {}Lessler, J.%
\end{APACrefauthors}%
\unskip\
\newblock
\APACrefYearMonthDay{2020}{}{}.
\newblock
{\BBOQ}\APACrefatitle {The incubation period of coronavirus disease 2019
  (COVID-19) from publicly reported confirmed cases: estimation and
  application} {The incubation period of coronavirus disease 2019 (covid-19)
  from publicly reported confirmed cases: estimation and application}.{\BBCQ}
\newblock
\APACjournalVolNumPages{Annals of internal medicine}{172}{9}{577--582}.
\PrintBackRefs{\CurrentBib}

\bibitem [\protect \citeauthoryear {%
Li%
\ \BBA {} Xiao%
}{%
Li%
\ \BBA {} Xiao%
}{%
{\protect \APACyear {2022}}%
}]{%
li2022complex}
\APACinsertmetastar {%
li2022complex}%
\begin{APACrefauthors}%
Li, T.%
\BCBT {}\ \BBA {} Xiao, Y.%
\end{APACrefauthors}%
\unskip\
\newblock
\APACrefYearMonthDay{2022}{}{}.
\newblock
{\BBOQ}\APACrefatitle {Complex dynamics of an epidemic model with saturated
  media coverage and recovery} {Complex dynamics of an epidemic model with
  saturated media coverage and recovery}.{\BBCQ}
\newblock
\APACjournalVolNumPages{Nonlinear Dynamics}{}{}{1--29}.
\PrintBackRefs{\CurrentBib}

\bibitem [\protect \citeauthoryear {%
Lim%
\ \protect \BOthers {.}}{%
Lim%
\ \protect \BOthers {.}}{%
{\protect \APACyear {2023}}%
}]{%
lim2023similar}
\APACinsertmetastar {%
lim2023similar}%
\begin{APACrefauthors}%
Lim, T\BPBI Y.%
, Xu, R.%
, Ruktanonchai, N.%
, Saucedo, O.%
, Childs, L\BPBI M.%
, Jalali, M\BPBI S.%
\BDBL {}Ghaffarzadegan, N.%
\end{APACrefauthors}%
\unskip\
\newblock
\APACrefYearMonthDay{2023}{}{}.
\newblock
{\BBOQ}\APACrefatitle {Why Similar Policies Resulted In Different COVID-19
  Outcomes: How Responsiveness And Culture Influenced Mortality Rates: Study
  examines why similar policies resulted in different COVID-19 outcomes in
  using data from more than 100 countries} {Why similar policies resulted in
  different covid-19 outcomes: How responsiveness and culture influenced
  mortality rates: Study examines why similar policies resulted in different
  covid-19 outcomes in using data from more than 100 countries}.{\BBCQ}
\newblock
\APACjournalVolNumPages{Health Affairs}{42}{12}{1637--1646}.
\PrintBackRefs{\CurrentBib}

\bibitem [\protect \citeauthoryear {%
Macdonald%
, Browne%
\BCBL {}\ \BBA {} Gulbudak%
}{%
Macdonald%
\ \protect \BOthers {.}}{%
{\protect \APACyear {2021}}%
}]{%
macdonald2021modelling}
\APACinsertmetastar {%
macdonald2021modelling}%
\begin{APACrefauthors}%
Macdonald, J.%
, Browne, C.%
\BCBL {}\ \BBA {} Gulbudak, H.%
\end{APACrefauthors}%
\unskip\
\newblock
\APACrefYearMonthDay{2021}{}{}.
\newblock
{\BBOQ}\APACrefatitle {Modelling {COVID-19 outbreaks in USA} with distinct
  testing, lockdown speed and fatigue rates} {Modelling {COVID-19 outbreaks in
  USA} with distinct testing, lockdown speed and fatigue rates}.{\BBCQ}
\newblock
\APACjournalVolNumPages{Royal Society Open Science}{8}{8}{210227}.
\PrintBackRefs{\CurrentBib}

\bibitem [\protect \citeauthoryear {%
Manrubia%
\ \BBA {} Zanette%
}{%
Manrubia%
\ \BBA {} Zanette%
}{%
{\protect \APACyear {2022}}%
}]{%
manrubia2022individual}
\APACinsertmetastar {%
manrubia2022individual}%
\begin{APACrefauthors}%
Manrubia, S.%
\BCBT {}\ \BBA {} Zanette, D\BPBI H.%
\end{APACrefauthors}%
\unskip\
\newblock
\APACrefYearMonthDay{2022}{}{}.
\newblock
{\BBOQ}\APACrefatitle {Individual risk-aversion responses tune epidemics to
  critical transmissibility {($R$= 1)}} {Individual risk-aversion responses
  tune epidemics to critical transmissibility {($R$= 1)}}.{\BBCQ}
\newblock
\APACjournalVolNumPages{Royal Society Open Science}{9}{4}{211667}.
\PrintBackRefs{\CurrentBib}

\bibitem [\protect \citeauthoryear {%
Martcheva%
}{%
Martcheva%
}{%
{\protect \APACyear {2015}}%
}]{%
martcheva2015introduction}
\APACinsertmetastar {%
martcheva2015introduction}%
\begin{APACrefauthors}%
Martcheva, M.%
\end{APACrefauthors}%
\unskip\
\newblock
\APACrefYear{2015}.
\newblock
\APACrefbtitle {An introduction to mathematical epidemiology} {An introduction
  to mathematical epidemiology}\ (\BVOL~61).
\newblock
\APACaddressPublisher{}{Springer}.
\PrintBackRefs{\CurrentBib}

\bibitem [\protect \citeauthoryear {%
Martcheva%
, Tuncer%
\BCBL {}\ \BBA {} Ngonghala%
}{%
Martcheva%
\ \protect \BOthers {.}}{%
{\protect \APACyear {2021}}%
}]{%
Martcheva}
\APACinsertmetastar {%
Martcheva}%
\begin{APACrefauthors}%
Martcheva, M.%
, Tuncer, N.%
\BCBL {}\ \BBA {} Ngonghala, C\BPBI N.%
\end{APACrefauthors}%
\unskip\
\newblock
\APACrefYearMonthDay{2021}{}{}.
\newblock
{\BBOQ}\APACrefatitle {Effects of social-distancing on infectious disease
  dynamics: an evolutionary game theory and economic perspective} {Effects of
  social-distancing on infectious disease dynamics: an evolutionary game theory
  and economic perspective}.{\BBCQ}
\newblock
\APACjournalVolNumPages{Journal of Biological Dynamics}{15}{1}{342--366}.
\PrintBackRefs{\CurrentBib}

\bibitem [\protect \citeauthoryear {%
Montefusco%
\ \protect \BOthers {.}}{%
Montefusco%
\ \protect \BOthers {.}}{%
{\protect \APACyear {2022}}%
}]{%
montefusco2022interacting}
\APACinsertmetastar {%
montefusco2022interacting}%
\begin{APACrefauthors}%
Montefusco, F.%
, Procopio, A.%
, Bulai, I\BPBI M.%
, Amato, F.%
, Pedersen, M\BPBI G.%
\BCBL {}\ \BBA {} Cosentino, C.%
\end{APACrefauthors}%
\unskip\
\newblock
\APACrefYearMonthDay{2022}{}{}.
\newblock
{\BBOQ}\APACrefatitle {Interacting with {COVID-19}: How population behavior,
  feedback and memory shaped recurrent waves of the epidemic} {Interacting with
  {COVID-19}: How population behavior, feedback and memory shaped recurrent
  waves of the epidemic}.{\BBCQ}
\newblock
\APACjournalVolNumPages{IEEE Control Systems Letters}{7}{}{583--588}.
\PrintBackRefs{\CurrentBib}

\bibitem [\protect \citeauthoryear {%
Morsky%
, Magpantay%
, Day%
\BCBL {}\ \BBA {} Ak{\c{c}}ay%
}{%
Morsky%
\ \protect \BOthers {.}}{%
{\protect \APACyear {2023}}%
}]{%
morsky2023impact}
\APACinsertmetastar {%
morsky2023impact}%
\begin{APACrefauthors}%
Morsky, B.%
, Magpantay, F.%
, Day, T.%
\BCBL {}\ \BBA {} Ak{\c{c}}ay, E.%
\end{APACrefauthors}%
\unskip\
\newblock
\APACrefYearMonthDay{2023}{}{}.
\newblock
{\BBOQ}\APACrefatitle {The impact of threshold decision mechanisms of
  collective behavior on disease spread} {The impact of threshold decision
  mechanisms of collective behavior on disease spread}.{\BBCQ}
\newblock
\APACjournalVolNumPages{Proceedings of the National Academy of
  Sciences}{120}{19}{e2221479120}.
\PrintBackRefs{\CurrentBib}

\bibitem [\protect \citeauthoryear {%
Mumtaz%
, Green%
\BCBL {}\ \BBA {} Duggan%
}{%
Mumtaz%
\ \protect \BOthers {.}}{%
{\protect \APACyear {2022}}%
}]{%
mumtaz2022exploring}
\APACinsertmetastar {%
mumtaz2022exploring}%
\begin{APACrefauthors}%
Mumtaz, N.%
, Green, C.%
\BCBL {}\ \BBA {} Duggan, J.%
\end{APACrefauthors}%
\unskip\
\newblock
\APACrefYearMonthDay{2022}{}{}.
\newblock
{\BBOQ}\APACrefatitle {Exploring the effect of misinformation on infectious
  disease transmission} {Exploring the effect of misinformation on infectious
  disease transmission}.{\BBCQ}
\newblock
\APACjournalVolNumPages{Systems}{10}{2}{50}.
\PrintBackRefs{\CurrentBib}

\bibitem [\protect \citeauthoryear {%
N'konzi%
, Chukwu%
\BCBL {}\ \BBA {} Nyabadza%
}{%
N'konzi%
\ \protect \BOthers {.}}{%
{\protect \APACyear {2022}}%
}]{%
n2022effect}
\APACinsertmetastar {%
n2022effect}%
\begin{APACrefauthors}%
N'konzi, J\BHBI P\BPBI N.%
, Chukwu, C\BPBI W.%
\BCBL {}\ \BBA {} Nyabadza, F.%
\end{APACrefauthors}%
\unskip\
\newblock
\APACrefYearMonthDay{2022}{}{}.
\newblock
{\BBOQ}\APACrefatitle {Effect of time-varying adherence to non-pharmaceutical
  interventions on the occurrence of multiple epidemic waves: A modeling study}
  {Effect of time-varying adherence to non-pharmaceutical interventions on the
  occurrence of multiple epidemic waves: A modeling study}.{\BBCQ}
\newblock
\APACjournalVolNumPages{Frontiers in Public Health}{10}{}{1087683}.
\PrintBackRefs{\CurrentBib}

\bibitem [\protect \citeauthoryear {%
Ochab%
, Manfredi%
, Puszynski%
\BCBL {}\ \BBA {} d’Onofrio%
}{%
Ochab%
\ \protect \BOthers {.}}{%
{\protect \APACyear {2023}}%
}]{%
ochab2023multiple}
\APACinsertmetastar {%
ochab2023multiple}%
\begin{APACrefauthors}%
Ochab, M.%
, Manfredi, P.%
, Puszynski, K.%
\BCBL {}\ \BBA {} d’Onofrio, A.%
\end{APACrefauthors}%
\unskip\
\newblock
\APACrefYearMonthDay{2023}{}{}.
\newblock
{\BBOQ}\APACrefatitle {Multiple epidemic waves as the outcome of stochastic
  {SIR} epidemics with behavioral responses: a hybrid modeling approach}
  {Multiple epidemic waves as the outcome of stochastic {SIR} epidemics with
  behavioral responses: a hybrid modeling approach}.{\BBCQ}
\newblock
\APACjournalVolNumPages{Nonlinear Dynamics}{111}{1}{887--926}.
\PrintBackRefs{\CurrentBib}

\bibitem [\protect \citeauthoryear {%
Pant%
, Safdar%
, Santillana%
\BCBL {}\ \BBA {} Gumel%
}{%
Pant%
\ \protect \BOthers {.}}{%
{\protect \APACyear {2024}}%
}]{%
pant2024mathematical}
\APACinsertmetastar {%
pant2024mathematical}%
\begin{APACrefauthors}%
Pant, B.%
, Safdar, S.%
, Santillana, M.%
\BCBL {}\ \BBA {} Gumel, A\BPBI B.%
\end{APACrefauthors}%
\unskip\
\newblock
\APACrefYearMonthDay{2024}{}{}.
\newblock
{\BBOQ}\APACrefatitle {{Mathematical Assessment of the Role of Human Behavior
  Changes on SARS-CoV-2 Transmission Dynamics in the United States}}
  {{Mathematical Assessment of the Role of Human Behavior Changes on SARS-CoV-2
  Transmission Dynamics in the United States}}.{\BBCQ}
\newblock
\APACjournalVolNumPages{Bulletin of Mathematical Biology}{86}{8}{92}.
\PrintBackRefs{\CurrentBib}

\bibitem [\protect \citeauthoryear {%
Perra%
}{%
Perra%
}{%
{\protect \APACyear {2021}}%
}]{%
perra2021non}
\APACinsertmetastar {%
perra2021non}%
\begin{APACrefauthors}%
Perra, N.%
\end{APACrefauthors}%
\unskip\
\newblock
\APACrefYearMonthDay{2021}{}{}.
\newblock
{\BBOQ}\APACrefatitle {Non-pharmaceutical interventions during the COVID-19
  pandemic: A review} {Non-pharmaceutical interventions during the covid-19
  pandemic: A review}.{\BBCQ}
\newblock
\APACjournalVolNumPages{Physics Reports}{913}{}{1--52}.
\PrintBackRefs{\CurrentBib}

\bibitem [\protect \citeauthoryear {%
Petropoulos%
\ \protect \BOthers {.}}{%
Petropoulos%
\ \protect \BOthers {.}}{%
{\protect \APACyear {2023}}%
}]{%
petropoulos2023operational}
\APACinsertmetastar {%
petropoulos2023operational}%
\begin{APACrefauthors}%
Petropoulos, F.%
, Laporte, G.%
, Aktas, E.%
, Alumur, S\BPBI A.%
, Archetti, C.%
, Ayhan, H.%
\BDBL {}others%
\end{APACrefauthors}%
\unskip\
\newblock
\APACrefYearMonthDay{2023}{}{}.
\newblock
{\BBOQ}\APACrefatitle {Operational Research: methods and applications}
  {Operational research: methods and applications}.{\BBCQ}
\newblock
\APACjournalVolNumPages{Journal of the Operational Research
  Society}{}{}{1--195}.
\PrintBackRefs{\CurrentBib}

\bibitem [\protect \citeauthoryear {%
Rahmandad%
, Lim%
\BCBL {}\ \BBA {} Sterman%
}{%
Rahmandad%
\ \protect \BOthers {.}}{%
{\protect \APACyear {2021}}%
}]{%
rahmandad2021behavioral}
\APACinsertmetastar {%
rahmandad2021behavioral}%
\begin{APACrefauthors}%
Rahmandad, H.%
, Lim, T\BPBI Y.%
\BCBL {}\ \BBA {} Sterman, J.%
\end{APACrefauthors}%
\unskip\
\newblock
\APACrefYearMonthDay{2021}{}{}.
\newblock
{\BBOQ}\APACrefatitle {Behavioral dynamics of {COVID-19}: estimating
  underreporting, multiple waves, and adherence fatigue across 92 nations}
  {Behavioral dynamics of {COVID-19}: estimating underreporting, multiple
  waves, and adherence fatigue across 92 nations}.{\BBCQ}
\newblock
\APACjournalVolNumPages{System Dynamics Review}{37}{1}{5--31}.
\PrintBackRefs{\CurrentBib}

\bibitem [\protect \citeauthoryear {%
Rahmandad%
\ \BBA {} Sterman%
}{%
Rahmandad%
\ \BBA {} Sterman%
}{%
{\protect \APACyear {2022}}%
}]{%
rahmandad2022quantifying}
\APACinsertmetastar {%
rahmandad2022quantifying}%
\begin{APACrefauthors}%
Rahmandad, H.%
\BCBT {}\ \BBA {} Sterman, J.%
\end{APACrefauthors}%
\unskip\
\newblock
\APACrefYearMonthDay{2022}{}{}.
\newblock
{\BBOQ}\APACrefatitle {Quantifying the {COVID-19} endgame: Is a new normal
  within reach?} {Quantifying the {COVID-19} endgame: Is a new normal within
  reach?}{\BBCQ}
\newblock
\APACjournalVolNumPages{System Dynamics Review}{38}{4}{329--353}.
\PrintBackRefs{\CurrentBib}

\bibitem [\protect \citeauthoryear {%
Rahmandad%
, Xu%
\BCBL {}\ \BBA {} Ghaffarzadegan%
}{%
Rahmandad%
\ \protect \BOthers {.}}{%
{\protect \APACyear {2022}}%
{\protect \APACexlab {{\protect \BCnt {1}}}}}]{%
rahmandad2022enhancing}
\APACinsertmetastar {%
rahmandad2022enhancing}%
\begin{APACrefauthors}%
Rahmandad, H.%
, Xu, R.%
\BCBL {}\ \BBA {} Ghaffarzadegan, N.%
\end{APACrefauthors}%
\unskip\
\newblock
\APACrefYearMonthDay{2022{\protect \BCnt {1}}}{}{}.
\newblock
{\BBOQ}\APACrefatitle {Enhancing long-term forecasting: Learning from {COVID-19
  models}} {Enhancing long-term forecasting: Learning from {COVID-19
  models}}.{\BBCQ}
\newblock
\APACjournalVolNumPages{PLOS Computational Biology}{18}{5}{e1010100}.
\PrintBackRefs{\CurrentBib}

\bibitem [\protect \citeauthoryear {%
Rahmandad%
, Xu%
\BCBL {}\ \BBA {} Ghaffarzadegan%
}{%
Rahmandad%
\ \protect \BOthers {.}}{%
{\protect \APACyear {2022}}%
{\protect \APACexlab {{\protect \BCnt {2}}}}}]{%
rahmandad2022missing}
\APACinsertmetastar {%
rahmandad2022missing}%
\begin{APACrefauthors}%
Rahmandad, H.%
, Xu, R.%
\BCBL {}\ \BBA {} Ghaffarzadegan, N.%
\end{APACrefauthors}%
\unskip\
\newblock
\APACrefYearMonthDay{2022{\protect \BCnt {2}}}{}{}.
\newblock
{\BBOQ}\APACrefatitle {A missing behavioural feedback in {COVID-19} models is
  the key to several puzzles} {A missing behavioural feedback in {COVID-19}
  models is the key to several puzzles}.{\BBCQ}
\newblock
\APACjournalVolNumPages{BMJ Global Health}{7}{10}{e010463}.
\PrintBackRefs{\CurrentBib}

\bibitem [\protect \citeauthoryear {%
Rai%
, Khajanchi%
, Tiwari%
, Venturino%
\BCBL {}\ \BBA {} Misra%
}{%
Rai%
\ \protect \BOthers {.}}{%
{\protect \APACyear {2022}}%
}]{%
rai2022impact}
\APACinsertmetastar {%
rai2022impact}%
\begin{APACrefauthors}%
Rai, R\BPBI K.%
, Khajanchi, S.%
, Tiwari, P\BPBI K.%
, Venturino, E.%
\BCBL {}\ \BBA {} Misra, A\BPBI K.%
\end{APACrefauthors}%
\unskip\
\newblock
\APACrefYearMonthDay{2022}{}{}.
\newblock
{\BBOQ}\APACrefatitle {Impact of social media advertisements on the
  transmission dynamics of {COVID-19 pandemic in India}} {Impact of social
  media advertisements on the transmission dynamics of {COVID-19 pandemic in
  India}}.{\BBCQ}
\newblock
\APACjournalVolNumPages{Journal of Applied Mathematics and
  Computing}{}{}{1--26}.
\PrintBackRefs{\CurrentBib}

\bibitem [\protect \citeauthoryear {%
Richardson%
}{%
Richardson%
}{%
{\protect \APACyear {2011}}%
}]{%
richardson2011reflections}
\APACinsertmetastar {%
richardson2011reflections}%
\begin{APACrefauthors}%
Richardson, G\BPBI P.%
\end{APACrefauthors}%
\unskip\
\newblock
\APACrefYearMonthDay{2011}{}{}.
\newblock
{\BBOQ}\APACrefatitle {Reflections on the foundations of system dynamics}
  {Reflections on the foundations of system dynamics}.{\BBCQ}
\newblock
\APACjournalVolNumPages{System dynamics review}{27}{3}{219--243}.
\PrintBackRefs{\CurrentBib}

\bibitem [\protect \citeauthoryear {%
Ross%
}{%
Ross%
}{%
{\protect \APACyear {1916}}%
}]{%
ross1916application}
\APACinsertmetastar {%
ross1916application}%
\begin{APACrefauthors}%
Ross, R.%
\end{APACrefauthors}%
\unskip\
\newblock
\APACrefYearMonthDay{1916}{}{}.
\newblock
{\BBOQ}\APACrefatitle {An application of the theory of probabilities to the
  study of a priori pathometry.—Part I} {An application of the theory of
  probabilities to the study of a priori pathometry.—part i}.{\BBCQ}
\newblock
\APACjournalVolNumPages{Proceedings of the Royal Society of London. Series A,
  Containing papers of a mathematical and physical
  character}{92}{638}{204--230}.
\PrintBackRefs{\CurrentBib}

\bibitem [\protect \citeauthoryear {%
Ross%
\ \BBA {} Hudson%
}{%
Ross%
\ \BBA {} Hudson%
}{%
{\protect \APACyear {1917}}%
{\protect \APACexlab {{\protect \BCnt {1}}}}}]{%
ross1917applicationII}
\APACinsertmetastar {%
ross1917applicationII}%
\begin{APACrefauthors}%
Ross, R.%
\BCBT {}\ \BBA {} Hudson, H\BPBI P.%
\end{APACrefauthors}%
\unskip\
\newblock
\APACrefYearMonthDay{1917{\protect \BCnt {1}}}{}{}.
\newblock
{\BBOQ}\APACrefatitle {An application of the theory of probabilities to the
  study of a priori pathometry.—Part II} {An application of the theory of
  probabilities to the study of a priori pathometry.—part ii}.{\BBCQ}
\newblock
\APACjournalVolNumPages{Proceedings of the Royal Society of London. Series A,
  Containing papers of a mathematical and physical
  character}{93}{650}{212--225}.
\PrintBackRefs{\CurrentBib}

\bibitem [\protect \citeauthoryear {%
Ross%
\ \BBA {} Hudson%
}{%
Ross%
\ \BBA {} Hudson%
}{%
{\protect \APACyear {1917}}%
{\protect \APACexlab {{\protect \BCnt {2}}}}}]{%
ross1917applicationIII}
\APACinsertmetastar {%
ross1917applicationIII}%
\begin{APACrefauthors}%
Ross, R.%
\BCBT {}\ \BBA {} Hudson, H\BPBI P.%
\end{APACrefauthors}%
\unskip\
\newblock
\APACrefYearMonthDay{1917{\protect \BCnt {2}}}{}{}.
\newblock
{\BBOQ}\APACrefatitle {An application of the theory of probabilities to the
  study of a priori pathometry.—Part III} {An application of the theory of
  probabilities to the study of a priori pathometry.—part iii}.{\BBCQ}
\newblock
\APACjournalVolNumPages{Proceedings of the Royal Society of London. Series A,
  Containing papers of a mathematical and physical
  character}{93}{650}{225--240}.
\PrintBackRefs{\CurrentBib}

\bibitem [\protect \citeauthoryear {%
Santos%
, Korah%
, Subramanian%
, Murugappan%
\BCBL {}\ \BBA {} Santos%
}{%
Santos%
\ \protect \BOthers {.}}{%
{\protect \APACyear {2021}}%
}]{%
santos2021infusing}
\APACinsertmetastar {%
santos2021infusing}%
\begin{APACrefauthors}%
Santos, E\BPBI E.%
, Korah, J.%
, Subramanian, S.%
, Murugappan, V.%
\BCBL {}\ \BBA {} Santos, E.%
\end{APACrefauthors}%
\unskip\
\newblock
\APACrefYearMonthDay{2021}{}{}.
\newblock
{\BBOQ}\APACrefatitle {Infusing culture in compartmental epidemic models}
  {Infusing culture in compartmental epidemic models}.{\BBCQ}
\newblock
\BIn{} \APACrefbtitle {2021 IEEE EMBS International Conference on Biomedical
  and Health Informatics (BHI)} {2021 ieee embs international conference on
  biomedical and health informatics (bhi)}\ (\BPGS\ 1--6).
\PrintBackRefs{\CurrentBib}

\bibitem [\protect \citeauthoryear {%
Shankar%
\ \protect \BOthers {.}}{%
Shankar%
\ \protect \BOthers {.}}{%
{\protect \APACyear {2021}}%
}]{%
shankar2021systematic}
\APACinsertmetastar {%
shankar2021systematic}%
\begin{APACrefauthors}%
Shankar, S.%
, Mohakuda, S\BPBI S.%
, Kumar, A.%
, Nazneen, P.%
, Yadav, A\BPBI K.%
, Chatterjee, K.%
\BCBL {}\ \BBA {} Chatterjee, K.%
\end{APACrefauthors}%
\unskip\
\newblock
\APACrefYearMonthDay{2021}{}{}.
\newblock
{\BBOQ}\APACrefatitle {Systematic review of predictive mathematical models of
  COVID-19 epidemic} {Systematic review of predictive mathematical models of
  covid-19 epidemic}.{\BBCQ}
\newblock
\APACjournalVolNumPages{Medical journal armed forces India}{77}{}{S385--S392}.
\PrintBackRefs{\CurrentBib}

\bibitem [\protect \citeauthoryear {%
Song%
\ \BBA {} Xiao%
}{%
Song%
\ \BBA {} Xiao%
}{%
{\protect \APACyear {2022}}%
}]{%
song2022analysis}
\APACinsertmetastar {%
song2022analysis}%
\begin{APACrefauthors}%
Song, P.%
\BCBT {}\ \BBA {} Xiao, Y.%
\end{APACrefauthors}%
\unskip\
\newblock
\APACrefYearMonthDay{2022}{}{}.
\newblock
{\BBOQ}\APACrefatitle {Analysis of a diffusive epidemic system with spatial
  heterogeneity and lag effect of media impact} {Analysis of a diffusive
  epidemic system with spatial heterogeneity and lag effect of media
  impact}.{\BBCQ}
\newblock
\APACjournalVolNumPages{Journal of Mathematical Biology}{85}{2}{17}.
\PrintBackRefs{\CurrentBib}

\bibitem [\protect \citeauthoryear {%
Sooknanan%
\ \BBA {} Seemungal%
}{%
Sooknanan%
\ \BBA {} Seemungal%
}{%
{\protect \APACyear {2023}}%
}]{%
sooknanan2023fomo}
\APACinsertmetastar {%
sooknanan2023fomo}%
\begin{APACrefauthors}%
Sooknanan, J.%
\BCBT {}\ \BBA {} Seemungal, T\BPBI A.%
\end{APACrefauthors}%
\unskip\
\newblock
\APACrefYearMonthDay{2023}{}{}.
\newblock
{\BBOQ}\APACrefatitle {FOMO (fate of online media only) in infectious disease
  modeling: a review of compartmental models} {Fomo (fate of online media only)
  in infectious disease modeling: a review of compartmental models}.{\BBCQ}
\newblock
\APACjournalVolNumPages{International Journal of Dynamics and
  Control}{11}{2}{892--899}.
\PrintBackRefs{\CurrentBib}

\bibitem [\protect \citeauthoryear {%
Sterman%
}{%
Sterman%
}{%
{\protect \APACyear {2000}}%
}]{%
sterman2000business}
\APACinsertmetastar {%
sterman2000business}%
\begin{APACrefauthors}%
Sterman, J\BPBI D.%
\end{APACrefauthors}%
\unskip\
\newblock
\APACrefYear{2000}.
\newblock
\APACrefbtitle {Business Dynamics: Systems Thinking and Modeling for a Complex
  World} {Business dynamics: Systems thinking and modeling for a complex
  world}.
\newblock
\APACaddressPublisher{}{McGraw-Hill}.
\PrintBackRefs{\CurrentBib}

\bibitem [\protect \citeauthoryear {%
Verelst%
, Willem%
\BCBL {}\ \BBA {} Beutels%
}{%
Verelst%
\ \protect \BOthers {.}}{%
{\protect \APACyear {2016}}%
}]{%
verelst2016behavioural}
\APACinsertmetastar {%
verelst2016behavioural}%
\begin{APACrefauthors}%
Verelst, F.%
, Willem, L.%
\BCBL {}\ \BBA {} Beutels, P.%
\end{APACrefauthors}%
\unskip\
\newblock
\APACrefYearMonthDay{2016}{}{}.
\newblock
{\BBOQ}\APACrefatitle {Behavioural change models for infectious disease
  transmission: a systematic review (2010--2015)} {Behavioural change models
  for infectious disease transmission: a systematic review
  (2010--2015)}.{\BBCQ}
\newblock
\APACjournalVolNumPages{Journal of The Royal Society
  Interface}{13}{125}{20160820}.
\PrintBackRefs{\CurrentBib}

\bibitem [\protect \citeauthoryear {%
Vynnycky%
\ \BBA {} White%
}{%
Vynnycky%
\ \BBA {} White%
}{%
{\protect \APACyear {2010}}%
}]{%
vynnycky2010introduction}
\APACinsertmetastar {%
vynnycky2010introduction}%
\begin{APACrefauthors}%
Vynnycky, E.%
\BCBT {}\ \BBA {} White, R.%
\end{APACrefauthors}%
\unskip\
\newblock
\APACrefYear{2010}.
\newblock
\APACrefbtitle {An introduction to infectious disease modelling} {An
  introduction to infectious disease modelling}.
\newblock
\APACaddressPublisher{}{Oxford University Press}.
\PrintBackRefs{\CurrentBib}

\bibitem [\protect \citeauthoryear {%
Wang%
, Zhang%
, Yan%
, Bai%
\BCBL {}\ \BBA {} He%
}{%
Wang%
\ \protect \BOthers {.}}{%
{\protect \APACyear {2023}}%
}]{%
wang2023evaluating}
\APACinsertmetastar {%
wang2023evaluating}%
\begin{APACrefauthors}%
Wang, A.%
, Zhang, X.%
, Yan, R.%
, Bai, D.%
\BCBL {}\ \BBA {} He, J.%
\end{APACrefauthors}%
\unskip\
\newblock
\APACrefYearMonthDay{2023}{}{}.
\newblock
{\BBOQ}\APACrefatitle {Evaluating the impact of multiple factors on the control
  of COVID-19 epidemic: {A modelling analysis using India as a case study}}
  {Evaluating the impact of multiple factors on the control of covid-19
  epidemic: {A modelling analysis using India as a case study}}.{\BBCQ}
\newblock
\APACjournalVolNumPages{Mathematical Biosciences and
  Engineering}{20}{4}{6237--6272}.
\PrintBackRefs{\CurrentBib}

\bibitem [\protect \citeauthoryear {%
Weitz%
, Park%
, Eksin%
\BCBL {}\ \BBA {} Dushoff%
}{%
Weitz%
\ \protect \BOthers {.}}{%
{\protect \APACyear {2020}}%
}]{%
weitz2020awareness}
\APACinsertmetastar {%
weitz2020awareness}%
\begin{APACrefauthors}%
Weitz, J\BPBI S.%
, Park, S\BPBI W.%
, Eksin, C.%
\BCBL {}\ \BBA {} Dushoff, J.%
\end{APACrefauthors}%
\unskip\
\newblock
\APACrefYearMonthDay{2020}{}{}.
\newblock
{\BBOQ}\APACrefatitle {Awareness-driven behavior changes can shift the shape of
  epidemics away from peaks and toward plateaus, shoulders, and oscillations}
  {Awareness-driven behavior changes can shift the shape of epidemics away from
  peaks and toward plateaus, shoulders, and oscillations}.{\BBCQ}
\newblock
\APACjournalVolNumPages{Proceedings of the National Academy of
  Sciences}{117}{51}{32764--32771}.
\PrintBackRefs{\CurrentBib}

\bibitem [\protect \citeauthoryear {%
Wirtz%
}{%
Wirtz%
}{%
{\protect \APACyear {2021}}%
}]{%
wirtz2021changing}
\APACinsertmetastar {%
wirtz2021changing}%
\begin{APACrefauthors}%
Wirtz, K.%
\end{APACrefauthors}%
\unskip\
\newblock
\APACrefYearMonthDay{2021}{}{}.
\newblock
{\BBOQ}\APACrefatitle {Changing readiness to mitigate {SARS}-{C}o{V}-2 steered
  long-term epidemic and social trajectories} {Changing readiness to mitigate
  {SARS}-{C}o{V}-2 steered long-term epidemic and social trajectories}.{\BBCQ}
\newblock
\APACjournalVolNumPages{Scientific Reports}{11}{1}{13919}.
\PrintBackRefs{\CurrentBib}

\bibitem [\protect \citeauthoryear {%
Yan%
, Cheke%
\BCBL {}\ \BBA {} Tang%
}{%
Yan%
\ \protect \BOthers {.}}{%
{\protect \APACyear {2023}}%
}]{%
yan2023coupling}
\APACinsertmetastar {%
yan2023coupling}%
\begin{APACrefauthors}%
Yan, Q.%
, Cheke, R\BPBI A.%
\BCBL {}\ \BBA {} Tang, S.%
\end{APACrefauthors}%
\unskip\
\newblock
\APACrefYearMonthDay{2023}{}{}.
\newblock
{\BBOQ}\APACrefatitle {{Coupling an individual adaptive-decision model with a
  SIRV model of influenza vaccination reveals new insights for epidemic
  control}} {{Coupling an individual adaptive-decision model with a SIRV model
  of influenza vaccination reveals new insights for epidemic control}}.{\BBCQ}
\newblock
\APACjournalVolNumPages{Statistics in Medicine}{42}{5}{716--729}.
\PrintBackRefs{\CurrentBib}

\bibitem [\protect \citeauthoryear {%
Zhang%
, Scarabel%
, Murty%
\BCBL {}\ \BBA {} Wu%
}{%
Zhang%
\ \protect \BOthers {.}}{%
{\protect \APACyear {2023}}%
}]{%
zhang2023renewal}
\APACinsertmetastar {%
zhang2023renewal}%
\begin{APACrefauthors}%
Zhang, X.%
, Scarabel, F.%
, Murty, K.%
\BCBL {}\ \BBA {} Wu, J.%
\end{APACrefauthors}%
\unskip\
\newblock
\APACrefYearMonthDay{2023}{}{}.
\newblock
{\BBOQ}\APACrefatitle {Renewal equations for delayed population behaviour
  adaptation coupled with disease transmission dynamics: A mechanism for
  multiple waves of emerging infections} {Renewal equations for delayed
  population behaviour adaptation coupled with disease transmission dynamics: A
  mechanism for multiple waves of emerging infections}.{\BBCQ}
\newblock
\APACjournalVolNumPages{Mathematical Biosciences}{365}{}{109068}.
\PrintBackRefs{\CurrentBib}

\bibitem [\protect \citeauthoryear {%
Zhou%
, Wang%
, Xia%
, Xiao%
\BCBL {}\ \BBA {} Tang%
}{%
Zhou%
\ \protect \BOthers {.}}{%
{\protect \APACyear {2020}}%
}]{%
zhou2020effects}
\APACinsertmetastar {%
zhou2020effects}%
\begin{APACrefauthors}%
Zhou, W.%
, Wang, A.%
, Xia, F.%
, Xiao, Y.%
\BCBL {}\ \BBA {} Tang, S.%
\end{APACrefauthors}%
\unskip\
\newblock
\APACrefYearMonthDay{2020}{}{}.
\newblock
{\BBOQ}\APACrefatitle {Effects of media reporting on mitigating spread of
  {COVID-19} in the early phase of the outbreak} {Effects of media reporting on
  mitigating spread of {COVID-19} in the early phase of the outbreak}.{\BBCQ}
\newblock
\APACjournalVolNumPages{Mathematical Biosciences and
  Engineering}{17}{3}{2693--2707}.
\PrintBackRefs{\CurrentBib}

\end{thebibliography}

\end{document}